\documentclass[a4paper,twoside,11pt]{amsart}

\usepackage[a4paper]{geometry}
\usepackage[latin1]{inputenc}
\usepackage[T1]{fontenc}

\usepackage{lipsum}
\usepackage{comment}
\usepackage{color}

\usepackage{amssymb,amsmath,amsthm,amscd,amsbsy}
\usepackage{mathrsfs}
\usepackage{stmaryrd}
\usepackage{esint}

\usepackage{enumitem}

\usepackage{subfig}

\usepackage{tabularx}
\usepackage{calc}
\usepackage[pdftex]{graphicx}
\usepackage[all]{xy}
\xyoption{v2}
\xyoption{2cell}
\UseAllTwocells

\usepackage[pdfauthor={Christophe Prange},pdftitle={Improved Regularity in Bumpy Lipschitz Domains},pdftex]{hyperref}

\numberwithin{equation}{section}

\DeclareMathOperator{\Idd}{I}

\DeclareMathOperator{\DtoN}{DN}

\DeclareMathOperator{\supp}{Supp}

\let\oldmarginpar\marginpar
\renewcommand\marginpar[1]{\oldmarginpar{\color{red}\raggedleft\tiny #1}}

\newcommand{\intbar}{{- \hspace{- 1.05 em}} \int} 

\title{Improved Regularity in Bumpy Lipschitz Domains}

\begin{document}

\newtheorem{theo}{Theorem}
\newtheorem{prop}[theo]{Proposition}
\newtheorem{lem}[theo]{Lemma}
\newtheorem{cor}[theo]{Corollary}
\newtheorem*{theo*}{Theorem}
\newtheorem{rst}{Result}
\renewcommand*{\therst}{\Alph{rst}}

\theoremstyle{definition}
\newtheorem{defi}[theo]{Definition}

\theoremstyle{remark}
\newtheorem{rem}{Remark}
\newtheorem*{rem*}{Remark}
\newtheorem*{rems*}{Remarks}

\title{Improved Regularity in Bumpy Lipschitz Domains}
\author{Carlos Kenig$^*$}
\thanks{$^*$The University of Chicago, $5734$ S. University Avenue, Chicago, IL $60637$, USA. \emph{E-mail address:} \texttt{cek@math.uchicago.edu}}
\author{Christophe Prange$^\dagger$}
\thanks{$^\dagger$The University of Chicago, $5734$ S. University Avenue, Chicago, IL $60637$, USA. \emph{E-mail address:} \texttt{cp@math.uchicago.edu}}

\begin{abstract}
This paper is devoted to the proof of Lipschitz regularity, down to the microscopic scale, for solutions of an elliptic system with highly oscillating coefficients, over a highly oscillating Lipschitz boundary. The originality of this result is that it does not assume more than Lipschitz regularity on the boundary. Our Theorem, which is a significant improvement of our previous work on Lipschitz estimates in bumpy domains, should be read as an improved regularity result for an elliptic system over a Lipschitz boundary. Our progress in this direction is made possible by an estimate for a boundary layer corrector. We believe that this estimate in the Sobolev-Kato class is of independent interest.
\end{abstract}

\maketitle

\pagestyle{plain}

\section{Introduction}

This paper is devoted to the proof of Lipschitz regularity, down to the microscopic scale, for weak solutions $u^\varepsilon=u^\varepsilon(x)\in\mathbb R^N$ of the elliptic system
\begin{equation}\label{sysuepsintro}
\left\{\begin{array}{rll}
-\nabla\cdot A(x/\varepsilon)\nabla u^\varepsilon&=0,&x\in D^\varepsilon_\psi(0,1),\\
u^\varepsilon&=0,&x\in\Delta^\varepsilon_\psi(0,1),
\end{array}\right.
\end{equation}
over a highly oscillating Lipschitz boundary. Throughout this work, $\psi$ is a Lipschitz graph,
\begin{equation*}
D^\varepsilon_{\psi}(0,1):=\{x'\in (-1,1)^{d-1},\ \varepsilon\psi(x'/\varepsilon)<x_d<\varepsilon\psi(x'/\varepsilon)+1\}\subset\mathbb R^d
\end{equation*}
and 
\begin{equation*}
\Delta^\varepsilon_{\psi}(0,1):=\{x'\in (-1,1)^{d-1},\ x_d=\varepsilon\psi(x'/\varepsilon)\}
\end{equation*}
is the lower highly oscillating boundary on which homogeneous Dirichlet boundary conditions are imposed. Our main theorem is the following.

\begin{theo}\label{theolipdowmicro}
There exists $C>0$ such that for all $\psi\in W^{1,\infty}(\mathbb R^{d-1})$, for all matrix $A=A(y)=(A^{\alpha\beta}_{ij}(y))\in\mathbb R^{d^2\times N^2}$, elliptic with constant $\lambda$, $1$-periodic and H\"older continuous with exponant $\nu>0$, for all $0<\varepsilon<1/2$, for all weak solution $u^\varepsilon$ to \eqref{sysuepsintro}, for all $r\in[\varepsilon,1/2]$
\begin{equation}\label{estlipdownmicroinpropintro}
\int_{(-r,r)^{d-1}}\int_{\varepsilon\psi(x'/\varepsilon)}^{\varepsilon\psi(x'/\varepsilon)+r}|\nabla u^\varepsilon|^2\leq Cr^d\int_{(-1,1)^{d-1}}\int_{\varepsilon\psi(x'/\varepsilon)}^{\varepsilon\psi(x'/\varepsilon)+1}|\nabla u^\varepsilon|^2,
\end{equation}
with $C=C(d,N,\lambda,[A]_{C^{0,\nu}},\|\psi\|_{W^{1,\infty}})$.
\end{theo}

The uniform estimate of Theorem \ref{theolipdowmicro} should be read as an improved regularity result. Indeed, estimate \eqref{estlipdownmicroinpropintro} can be seen as a Lipschitz estimate down to the microscopic scale $O(\varepsilon)$. Its originality lies in the fact that no smoothness of the boundary, which is just assumed to be Lipschitz, is needed for it to hold. Previous results in this direction always relied on some smoothness of the boundary, typically $\psi\in C^{1,\nu}$ with $\nu>0$, or $\psi\in C^1_\omega$ with $\omega$ a modulus of continuity satisfying a Dini type condition, i.e. $\int_0^1\omega(t)/tdt<\infty$. 

Pioneering work on uniform estimates in homogenization has been achieved by Avellaneda and Lin in the late 80's \cite{alin,alinscal,alin2,Alin90P,alinLp}. The regularity theory for operators with highly oscillating coefficients has recently attracted a lot of attention, and important contributions have been made to relax the structure assumptions on the oscillations \cite{ArmstrongShen14,ArSmart14,GloriaNeukammOtto14}. Our work is in a different vein. It is focused on the boundary behavior of solutions.

Theorem \ref{theolipdowmicro} represents a considerable improvement of a recent result obtained by the two authors, namely Result B and Theorem 16 in \cite{BLtailosc}. This first work dealt with uniform Lipschitz regularity over highly oscillating $C^{1,\nu}$ boundaries. To the best of our knowledge, an improved regularity result up to the boundary such as the one of Theorem \ref{theolipdowmicro} is new. Our breakthrough is made possible by estimating a boundary layer corrector $v=v(y)$ solution to the system
\begin{equation}\label{sysblbumpyhpintro}
\left\{
\begin{array}{rll}
-\nabla\cdot A(y)\nabla v&=0,&y_d>\psi(y'),\\
v&=v_0,&y_d=\psi(y'),
\end{array}
\right.
\end{equation}
in the Lipschitz half-space $y_d>\psi(y')$ with non localized Dirichlet boundary data $v_0$.

\begin{theo}\label{theoblbumpy}
Assume $\psi\in W^{1,\infty}(\mathbb R^{d-1})$ and $v_0\in H^{1/2}_{uloc}(\mathbb R^{d-1})$ i.e.
\begin{equation*}
\sup_{\xi\in\mathbb Z^{d-1}}\|v_0\|_{H^{1/2}(\xi+(0,1)^{d-1})}<\infty.
\end{equation*}
Then, there exists a unique weak solution $v$ of \eqref{sysblbumpyhpintro} such that 
\begin{equation}\label{estbumpyhpuloc}
\sup_{\xi\in\mathbb Z^{d-1}}\int_{\xi+(0,1)^{d-1}}\int_{\psi(y')}^\infty|\nabla v|^2dy_ddy'\leq C\|v_0\|_{H^{1/2}_{uloc}}^2<\infty,
\end{equation}
with $C=C(d,N,\lambda,[A]_{C^{0,\nu}},\|\psi\|_{W^{1,\infty}})$.
\end{theo}

\subsection*{Overview of the paper}

In section \ref{secprelim} we recall several results related to Sobolev-Kato spaces, homogenization and uniform Lipschitz estimates. These results are of constant use in our work. Then the paper has two main parts. The first aim is to prove Theorem \ref{theoblbumpy} about the well-posedness of the boundary layer system in a space of non localized energy over a Lipschitz boundary. The key idea is to carry out a domain decomposition. Subsequently, there are three steps. Firstly, we prove the well-posedness of the boundary layer system over a flat boundary, namely in the domain $\mathbb R^d_+$. This is done in section \ref{secflathalfspace}. Secondly, we define and estimate a Dirichlet to Neumann operator over $H^{1/2}_{uloc}$. This key tool is introduced in section \ref{secDN}. Thirdly, we show that proving the well-posedness of the boundary layer system over a Lipschitz boundary boils down to analyzing a problem in a layer $\{\psi(y')<y_d<0\}$ close to the boundary. The energy estimates for this problem are carried out in section \ref{secwpbumpyhp}. Eventually in section \ref{secimpro}, and this is the last part of this work, we are able to prove Theorem \ref{theolipdowmicro} using a compactness scheme.

\subsection*{Framework and notations}

Let $\lambda>0$ and $0<\nu<1$ be fixed in what follows. We assume that the coefficients matrix $A=A(y)=(A^{\alpha\beta}_{ij}(y))$, with $1\leq\alpha,\ \beta\leq d$ and $1\leq i,\ j\leq N$ is real, that
\begin{equation}\label{smoothA2}
A\in C^{0,\nu}(\mathbb R^{d}),
\end{equation}
that $A$ is uniformly elliptic i.e.
\begin{equation}\label{elliptA}
\lambda |\xi|^2\leq A^{\alpha\beta}_{ij}(y)\xi^\alpha_i\xi^\beta_j\leq \frac{1}{\lambda}|\xi|^2,\quad\mbox{for all}\ \xi=(\xi^\alpha_i)\in\mathbb R^{dN},\ y\in\mathbb R^d
\end{equation}
and periodic i.e. 
\begin{equation}\label{perA}
A(y+z)=A(y),\quad\mbox{for all}\ y\in\mathbb R^d,\ z\in\mathbb Z^d.
\end{equation}
We say that $A$ belongs to the class $\mathcal A^\nu$ if $A$ satisfies \eqref{smoothA2}, \eqref{elliptA} and \eqref{perA}.

For easy reference, we summarize here the standard notations used throughout the text. For $x\in\mathbb R^d$, $x=(x',x_d)$, so that $x'\in\mathbb R^{d-1}$ denotes the $d-1$ first components of the vector $x$. For $\varepsilon>0$, $r>0$, let
\begin{equation*}
\begin{aligned}
D^\varepsilon_{\psi}(0,r)&:=\left\{(x',x_d),\ |x'|<r,\ \varepsilon\psi(x'/\varepsilon)<x_d<\varepsilon\psi(x'/\varepsilon)+r\right\},\\
\Delta^\varepsilon_{\psi}(0,r)&:=\left\{(x',x_d),\ |x'|<r,\ x_d=\varepsilon\psi(x'/\varepsilon)\right\},\\
D_0(0,r)&:=\left\{(x',x_d),\ |x'|<r,\ 0<x_d<r\right\},\qquad \Delta_0(0,r):=\left\{(x',0),\ |x'|<r\right\},\\
\mathbb R^d_+&:=\mathbb R^{d-1}\times(0,\infty),\qquad\Omega_+:=\{\psi(y')<y_d\},\\
\Omega_\flat&:=\{\psi(y')<y_d<0\},\qquad\Sigma_k:=(-k,k)^{d-1},
\end{aligned}
\end{equation*}
where $|x'|=\max_{i=1,\ldots\ d}|x_i|$. We sometimes write $D_\psi(0,r)$ and $\Delta_\psi(0,r)$ in short for $D^1_\psi(0,r)$ and $\Delta^1_\psi(0,r)$; in that situation the boundary is not highly oscillating because $\varepsilon=1$. Let also
\begin{equation*}
(\overline{u})_{D^\varepsilon_\psi(0,r)}:=\intbar_{D^\varepsilon_\psi(0,r)}u=\frac{1}{|D^\varepsilon_\psi(0,r)|}\int_{D^\varepsilon_\psi(0,r)}u.
\end{equation*}
The Lebesgue measure of a set is denoted by $|\cdot|$. For a positive integer $m$, let also $\Idd_m$ denote the identity matrix $M_m(\mathbb R)$. The function $\mathbf{1}_E$ denotes the characteristic function of a set $E$. The notation $\eta$ usually stands for a cut-off function. Ad hoc definitions are given when needed. Unless stated otherwise, the duality product $\langle \cdot,\cdot\rangle:=\langle \cdot,\cdot\rangle_{\mathcal D',\mathcal D}$ always denotes the duality between $\mathcal D(\mathbb R^{d-1})=C^\infty_0(\mathbb R^{d-1})$ and $\mathcal D'$. In the sequel, $C>0$ is always a constant uniform in $\varepsilon$ which may change from line to line. 

\section{Preliminaries}
\label{secprelim}

\subsection{On Sobolev-Kato spaces} 

For $s\geq 0$, we define the Sobolev-Kato space $H^s_{uloc}(\mathbb R^{d-1})$ of functions of non localized $H^s$ energy by
\begin{equation*}
H^s_{uloc}(\mathbb R^{d-1}):=\left\{u\in H^s_{loc}(\mathbb R^d),\ \sup_{\xi\in \mathbb Z^{d-1}}\|u\|_{H^s(\xi+(0,1)^{d-1})}<\infty\right\}.
\end{equation*}
We will mainly work with $H^{1/2}_{uloc}$. The following lemma is a useful tool to compare the $H^{1/2}_{uloc}$ norm to the $H^{1/2}$ norm of a $H^{1/2}(\mathbb R^{d-1})$ function.

\begin{lem}\label{lemh1/2H1/2uloc}
Let $\eta\in C^\infty_c(\mathbb R^{d-1})$ and $v_0\in H^{1/2}_{uloc}(\mathbb R^{d-1})$. Assume that $\supp\eta\subset B(0,R)$, for $R>0$. Then,
\begin{equation}\label{estulocH12}
\|\eta v_0\|_{H^{1/2}}\leq CR^{\frac{d-1}{2}}\|v_0\|_{H^{1/2}_{uloc}},
\end{equation} 
with $C=C(d,\|\eta\|_{W^{1,\infty}})$.
\end{lem}

For a proof, we refer to the proof of Lemma 2.26 in \cite{DP_SC}.

\subsection{Homogenization and weak convergence}

We recall the standard weak convergence result in periodic homogenization for a fixed domain $\Omega$. As usual, the constant homogenized matrix $\overline{A}=\overline{A}^{\alpha\beta}\in M_N(\mathbb R)$ is given by
\begin{equation}\label{defA0}
\overline{A}^{\alpha\beta}:=\int_{\mathbb T^d}A^{\alpha\beta}(y)dy+\int_{\mathbb T^d}A^{\alpha\gamma}(y)\partial_{y_\gamma}\chi^\beta(y)dy,
\end{equation}
where the family $\chi=\chi^\gamma(y)\in M_N(\mathbb R)$, $y\in\mathbb T^d$, solves the cell problems
\begin{equation}\label{eqdefchi}
-\nabla_y \cdot A(y)\nabla_y \chi^\gamma=\partial_{y_\alpha}A^{\alpha\gamma},\ y\in \mathbb T^d\qquad \mbox{and}\qquad \int_{\mathbb T^d}\chi^\gamma(y)dy=0.
\end{equation}

\begin{theo}[weak convergence]\label{theoweakcvhomo}
Let $\Omega$ be a bounded Lipschitz domain in $\mathbb R^d$ and let $u_k\in H^1(\Omega)$ be a sequence of weak solutions to
\begin{equation*}
-\nabla\cdot A_k(x/\varepsilon_k)\nabla u_k=f_k\in (H^1(\Omega))',
\end{equation*}
where $\varepsilon_k\rightarrow 0$ and the matrices $A_k=A_k(y)\in L^\infty$ satisfy \eqref{elliptA} and \eqref{perA}. Assume that there exist $f\in (H^1(\Omega))'$ and $u_k\in W^{1,2}(\Omega)$, such that $f_k\longrightarrow f$ strongly in $(H^1(\Omega))'$, $u_k\rightarrow u_0$ strongly in $L^2(\Omega)$ and $\nabla u_k\rightharpoonup\nabla u^0$ weakly in $L^2(\Omega)$. Also assume that the constant matrix $\overline{A_k}$ defined by \eqref{defA0} with $A$ replaced by $A_k$ converges to a constant matrix $A^0$. Then
\begin{equation*}
A_k(x/\varepsilon_k)\nabla u_k\rightharpoonup A^0\nabla u^0\quad\mbox{weakly in}\ L^2(\Omega) 
\end{equation*}
and 
\begin{equation*}
\nabla\cdot A^0\nabla u^0=f\in (H^1(\Omega))'.
\end{equation*}
\end{theo}

For a proof, which relies on the classical oscillating test function argument, we refer for instance to \cite[Lemma 2.1]{KLSNeumann}. This is an interior convergence result, since no boundary condition is prescribed on $u_k$.

\subsection{Uniform estimates in homogenization and applications}

We recall here the boundary Lipschitz estimate proved by Avellaneda and Lin in \cite{alin}.

\begin{theo}[Lipschitz estimate, {\cite[Lemma 20]{alin}}]\label{theoestlipbdaryalin}
For all $\kappa>0$, $0<\mu<1$, there exists $C>0$ such that for all $\psi\in C^{1,\nu}(\mathbb R^{d-1})\cap W^{1,\infty}(\mathbb R^{d-1})$, for all $A\in\mathcal A^{\nu}$, for all $r>0$, for all $\varepsilon>0$, for all $f\in L^{d+\kappa}(D_\psi(0,r))$, for all $F\in C^{0,\mu}(D_\psi(0,r))$, for all $u^\varepsilon\in L^\infty(D_\psi(0,r))$ weak solution to
\begin{equation*}
\left\{\begin{array}{rll}
-\nabla\cdot A(x/\varepsilon)\nabla u^\varepsilon&=f+\nabla\cdot F&x\in D_{\psi}(0,r),\\
u^\varepsilon&=0,&x\in\Delta_\psi(0,r),
\end{array}\right.
\end{equation*}
the following estimate holds
\begin{equation}\label{estlipbdaryalin}
\|\nabla u^\varepsilon\|_{L^\infty(D_\psi(0,r/2))}\leq C\left\{r^{-1}\|u^\varepsilon\|_{L^\infty(D_\psi(0,r))}+r^{1-d/(d+\kappa)}\|f\|_{L^{d+\kappa}(D_\psi(0,r))}+r^{\mu}\|F\|_{C^{0,\mu}(D_\psi(0,r))}\right\}.
\end{equation}
Notice that $C=C(d,N,\lambda,\kappa,\mu,\|\psi\|_{W^{1,\infty}},[\nabla\psi]_{C^{0,\nu}},[A]_{C^{0,\nu}})$.
\end{theo}

As stated in our earlier work \cite{BLtailosc}, this estimate does not cover the case of highly oscillating boundaries, since the constant in \eqref{estlipbdaryalin} involves the $C^{0,\nu}$ semi-norm of $\nabla\psi$.

In this work, we rely on Theorem \ref{theoestlipbdaryalin} to get large-scale pointwise estimates on the Poisson kernel $P=P(y,\tilde{y})$ associated to the domain $\mathbb R^d_+$ and to the operator $-\nabla\cdot A(y)\nabla$.

\begin{prop}\label{propestpeps}
For all $d\geq 2$, there exists $C>0$, such that for all $A\in \mathcal A^{\nu}$, we have:
\begin{enumerate}[label=(\arabic*)]
\item for all $y\in \mathbb R^d_+$, for all $\tilde{y}\in \mathbb R^{d-1}\times\{0\}$, we have
\begin{align}
|P(y,\tilde{y})|&\leq\frac{Cy_d}{|y-\tilde{y}|^d},\label{estPoissonkernel}\\
|\nabla_yP(y,\tilde{y})|&\leq\frac{C}{|y-\tilde{y}|^d},\label{estnablaP}
\end{align}
\item for all $y,\ \tilde{y}\in \mathbb R^{d-1}\times\{0\}$, $y\neq\tilde{y}$,
\begin{equation}\label{estnablaPbdary}
|\nabla_yP(y,\tilde{y})|\leq\frac{C}{|y-\tilde{y}|^d}.
\end{equation}
\end{enumerate}
Notice that $C=C(d,N,\lambda,[A]_{C^{0,\nu}})$.
\end{prop}

The proof of those estimates starting from the uniform Lipschitz estimate of Theorem \ref{theoestlipbdaryalin} is standard (see for instance \cite{alin}).

\section{Boundary layer corrector in a flat half-space}
\label{secflathalfspace}

This section is devoted to the well-posedness of the boundary layer problem
\begin{equation}\label{sysblhp}
\left\{
\begin{array}{rll}
-\nabla\cdot A(y)\nabla v&=0,&y_d>0,\\
v&=v_0\in H^{1/2}_{uloc}(\mathbb R^{d-1}),&y_d=0,
\end{array}
\right.
\end{equation}
in the flat half-space $\mathbb R^d_+$.

\begin{theo}\label{theoblflat}
Assume $v_0\in H^{1/2}_{uloc}(\mathbb R^{d-1})$. Then, there exists a unique weak solution $v$ of \eqref{sysblbumpyhp} such that 
\begin{equation}\label{estulocflathp}
\sup_{\xi\in\mathbb Z^{d-1}}\int_{\xi+(0,1)^{d-1}}\int_{0}^\infty|\nabla v|^2dy_ddy'\leq C\|v_0\|_{H^{1/2}_{uloc}}^2<\infty,
\end{equation}
with $C=C(d,N,\lambda,[A]_{C^{0,\nu}})$.
\end{theo}

The proof is in three steps: (i) we define a function $v$ and prove it is a weak solution to \eqref{sysblhp}, (ii) we prove that the solution we have defined satisfies the estimate \eqref{estulocflathp}, (iii) we prove uniqueness of solutions verifying \eqref{estulocflathp}.

\subsection{Existence of a weak solution}

Let $\eta\in C^\infty_c(\mathbb R)$ a cut-off function such that 
\begin{equation}\label{defeta}
\eta\equiv 1\ \mbox{on}\ (-1,1),\qquad 0\leq\eta\leq 1,\qquad \|\eta'\|_{L^\infty}\leq 2.
\end{equation}
Let $y_*\in\mathbb R^d_+$ be fixed. Notice that
\begin{equation*}
\eta(|\cdot-y_*'|)\in C^\infty_c(\mathbb R^{d-1}),\quad \eta(|\cdot-y_*'|)\equiv 1\ \mbox{on}\ B(y_*',1),\quad 0\leq \eta(|\cdot-y_*'|)\leq 1\quad\mbox{and}\quad\|\nabla(\eta(|\cdot-y_*'|))\|_{L^\infty}\leq 2.
\end{equation*}
We define
\begin{equation}\label{defsolflateta}
v(y_*):=v^\sharp(y_*)+v^\flat(y_*),
\end{equation}
where for $y\in\mathbb R^d_+$,
\begin{equation*}
v^\sharp(y):=\int_{\mathbb R^{d-1}\times\{0\}}P(y,\tilde{y})(1-\eta(|\tilde{y}'-y_*'|))v_0(\tilde{y}')d\tilde{y},
\end{equation*}
and $v^\flat=v^\flat(y)\in H^1(\mathbb R^d_+)$ is the unique weak solution to 
\begin{equation*}
\left\{
\begin{array}{rll}
-\nabla\cdot A(y)\nabla v^\flat&=0,&y_d>0,\\
v^\flat&=\eta(|y'-y_*'|)v_0(y')\in H^{1/2}(\mathbb R^{d-1}),&y_d=0,
\end{array}
\right.
\end{equation*}
satisfying 
\begin{equation}\label{estvflatH1}
\int_{\mathbb R^{d}_+}|\nabla v^\flat|^2dy'dy_d\leq C\|\eta v_0\|_{H^{1/2}}^2,
\end{equation}
with $C=C(d,N,\lambda)$. First of all, one has to prove that the definition of $v$ does not depend on the choice of the cut-off $\eta$. Let $\eta_1,\ \eta_2\in C^\infty_c(\mathbb R)$ be two cut-off functions satisfying \eqref{defeta}. We denote by $v_1(y_*)$ and $v_2(y_*)$ the associated vectors defined by 
\begin{equation*}
\begin{aligned}
v_1(y_*)&:=\int_{\mathbb R^{d-1}\times\{0\}}P(y_*,\tilde{y})(1-\eta_1(|\tilde{y}'-y_*'|)v_0(\tilde{y}')d\tilde{y}+v^\flat_1(y_*),\\
v_2(y_*)&:=\int_{\mathbb R^{d-1}\times\{0\}}P(y_*,\tilde{y})(1-\eta_2(|\tilde{y}'-y_*'|))v_0(\tilde{y}')d\tilde{y}+v^\flat_2(y_*).
\end{aligned}
\end{equation*}
Substracting, we get
\begin{equation}\label{substrindepeta}
v_1(y_*)-v_2(y_*)=\int_{\mathbb R^{d-1}\times\{0\}}P(y_*,\tilde{y})(\eta_2(|\tilde{y}'-y_*'|)-\eta_1(|\tilde{y}'-y_*'|))v_0(\tilde{y}')d\tilde{y}+v^\flat_1(y_*)-v^\flat_2(y_*).
\end{equation}
Now since
\begin{equation*}
y\longmapsto\int_{\mathbb R^{d-1}\times\{0\}}P(y,\tilde{y})(\eta_2(|\tilde{y}'-y_*'|)-\eta_1(|\tilde{y}'-y_*'|))v_0(\tilde{y}')d\tilde{y}
\end{equation*}
is the unique solution to 
\begin{equation*}
\left\{
\begin{array}{rll}
-\nabla\cdot A(y)\nabla v^\flat&=0,&y_d>0,\\
v^\flat&=(\eta_2(|y'-y_*'|)-\eta_1(|y'-y_*'|))v_0(y')\in H^{1/2}(\mathbb R^{d-1}),&y_d=0,
\end{array}
\right.
\end{equation*}
the difference in \eqref{substrindepeta} has to be zero, which proves that our definition of $v$ is independent of the choice of $\eta$.

It remains to prove that $v=v(y)$ defined by \eqref{defsolflateta} is actually a weak solution to \eqref{sysblhp}. Let $\varphi_\diamond=\varphi_\diamond(y')\in C^\infty_c(\mathbb R^{d-1})$ and $\varphi_d=\varphi_d(y_d)\in C^\infty_c((0,\infty))$. We choose $\eta\in C^\infty_c(\mathbb R)$ satisfying \eqref{defeta} and such that $\eta(|\cdot|)\equiv 1$ on $\supp\varphi_\diamond+B(0,1)$. We aim at proving
\begin{equation*}
\int_{\mathbb R^{d}_+}v(y)\left(-\nabla\cdot A^*(y)\nabla(\varphi_\diamond\varphi_d)\right)dy=0.
\end{equation*}
This relation is clear for $v^\flat$. For $v^\sharp$, by Fubini and then integration by parts
\begin{equation*}
\begin{aligned}
&\int_{\mathbb R^d_+}v^\sharp(y)\left(-\nabla\cdot A^*(y)\nabla(\varphi_\diamond\varphi_d)\right)dy\\
&=\int_{\supp\varphi_\diamond\times\supp\varphi_d}\int_{\mathbb R^{d-1}\times\{0\}}P(y,\tilde{y})(1-\eta(\tilde{y}))v_0(\tilde{y}')\left(-\nabla\cdot A^*(y)\nabla\varphi_\diamond\varphi_d\right)d\tilde{y}dy\\
&=\int_{\mathbb R^{d-1}\times\{0\}}\int_{\supp\varphi_\diamond\times\supp\varphi_d}P(y,\tilde{y})\left(-\nabla\cdot A^*(y)\nabla(\varphi_\diamond\varphi_d)\right)dy(1-\eta(\tilde{y}))v_0(\tilde{y}')d\tilde{y}\\
&=\int_{\mathbb R^{d-1}\times\{0\}}\left\langle-\nabla\cdot A(y)\nabla P(y,\tilde{y}),\varphi_\diamond\varphi_d\right\rangle (1-\eta(\tilde{y}))v_0(\tilde{y}')d\tilde{y}=0.
\end{aligned}
\end{equation*}

\subsection{Gradient estimate}

Let $\varphi_\diamond=\varphi_\diamond(y')\in C^\infty_c(\mathbb R^{d-1})$ and $\varphi_d=\varphi_d(y_d)\in C^\infty_c((0,\infty))$. We choose $R>1$ such that $\supp\varphi_\diamond+B(0,1)\subset B(0,R)$. Our goal is to prove
\begin{equation*}
\left|\int_{\mathbb R^d_+}\nabla v(y)\varphi_\diamond\varphi_d(y)dy\right|\leq CR^{\frac{d-1}{2}}\|v_0\|_{H^{1/2}_{uloc}}\|\varphi_\diamond\|_{L^2}\|\varphi_d\|_{L^2},
\end{equation*}
with $C=C(d,N,\lambda,[A]_{C^{0,\nu}})$. This estimate clearly implies the bound \eqref{estulocflathp}. Let $\eta\in C^{\infty}_c(\mathbb R)$ such that \eqref{defeta}
\begin{equation*}
\eta(|\cdot|)\equiv 1\ \mbox{on}\  B(0,R)\quad\mbox{and}\quad\supp\eta(|\cdot|)\subset B(0,2R).
\end{equation*}
Combining \eqref{estvflatH1} and the result of Lemma \ref{lemh1/2H1/2uloc}, we get
\begin{equation*}
\int_{\mathbb R^{d}_+}|\nabla v^\flat|^2dy'dy_d\leq CR^{d-1}\|v_0\|^2_{H^{1/2}_{uloc}},
\end{equation*}
with $C=C(d,N,\lambda)$. 

It remains to estimate
\begin{equation*}
\int_{\mathbb R^d_+}\nabla v^\sharp(y)\varphi_\diamond\varphi_d(y)dy=\int_0^1\int_{\mathbb R^{d-1}}\nabla v^\sharp(y)\varphi_\diamond\varphi_d(y)dy'dy_d+\int_1^\infty\int_{\mathbb R^{d-1}}\nabla v^\sharp(y)\varphi_\diamond\varphi_d(y)dy'dy_d.
\end{equation*}
To estimate these terms we rely on the the bound \eqref{estnablaP}: for all $y\in\mathbb R^d_+$, $\tilde{y}\in\mathbb R^{d-1}\times\{0\}$,
\begin{equation*}
|\nabla_yP(y,\tilde{y})|\leq \frac{C}{|y-\tilde{y}|^d}=\frac{C}{(y_d^2+|y'-\tilde{y}'|^2)^{d/2}},
\end{equation*}
with $C=C(d,N,\lambda,[A]_{C^{0,\nu}})$.

We begin with two useful estimates. For $y\in B(0,R)$, we have on the one hand
\begin{equation}\label{usefulest1}
\begin{aligned}
\int_{\mathbb R^{d-1}}\frac{1}{|y'-\tilde{y}'|^d}(1-\eta(|\tilde{y}'|))|v_0(\tilde{y}')|^2d\tilde{y}&\leq \int_{\mathbb R^{d-1}\setminus B(0,1)}\frac{1}{|y'-\tilde{y}'|^d}|v_0(\tilde{y}')|^2d\tilde{y}\\
&\leq\sum_{\xi\in\mathbb Z^{d-1}\setminus\{0\}}\frac{1}{|\xi|^d}\|v_0\|_{L^2_{uloc}}^2\leq C\|v_0\|_{L^2_{uloc}}^2
\end{aligned}
\end{equation}
and on the other hand
\begin{equation}\label{usefulest2}
\begin{aligned}
\int_{\mathbb R^{d-1}}\frac{1}{(y_d^2+|y'-\tilde{y}'|^2)^{d/2}}(1-\eta(|\tilde{y}'|))|v_0(\tilde{y}')|^2d\tilde{y}&\leq \int_{\mathbb R^{d-1}}\frac{1}{(y_d^2+|y'-\tilde{y}'|^2)^{d/2}}|v_0(\tilde{y}')|^2d\tilde{y}\\
&\leq\int_{\mathbb R^{d-1}}\frac{1}{(y_d^2+|y'-\tilde{y}'|^2)^{d/2}}d\tilde{y}\|v_0\|_{L^2_{uloc}}^2\\
&\leq \frac{C}{y_d}\|v_0\|_{L^2_{uloc}}^2.
\end{aligned}
\end{equation}
Using \eqref{usefulest1}, we get
\begin{equation*}
\begin{aligned}
&\left|\int_0^1\int_{\mathbb R^{d-1}}\nabla v^\sharp(y)\varphi_\diamond\varphi_d(y)dy'dy_d\right|=\left|\int_0^1\int_{\mathbb R^{d-1}}\int_{\mathbb R^{d-1}\times\{0\}}\nabla_yP(y,\tilde{y})(1-\eta(|\tilde{y}'|))v_0(\tilde{y}')d\tilde{y}\varphi_\diamond\varphi_d(y)dy'dy_d\right|\\
&\leq \int_0^1\int_{\mathbb R^{d-1}}\left(\int_{\mathbb R^{d-1}\times\{0\}}\frac{1-\eta(|\tilde{y}'|)}{|y'-\tilde{y}'|^d}d\tilde{y}'\right)^{1/2}\left(\int_{\mathbb R^{d-1}\times\{0\}}\frac{1-\eta(|\tilde{y}'|)}{|y'-\tilde{y}'|^d}|v_0(\tilde{y}')|^2d\tilde{y}\right)^{1/2}|\varphi_\diamond\varphi_d(y)|dy'dy_d\\
&\leq C\|v_0\|_{L^2_{uloc}}\int_0^1\int_{\mathbb R^{d-1}}\left(\int_1^\infty\frac{1}{r^2}\right)^{1/2}|\varphi_\diamond\varphi_d(y)|dy'dy_d\\
&\leq C\|v_0\|_{L^2_{uloc}}\int_0^1\int_{\mathbb R^{d-1}}|\varphi_\diamond\varphi_d(y)|dy'dy_d\leq CR^{\frac{d-1}{2}}\|v_0\|_{L^2_{uloc}}\|\varphi_\diamond\|_{L^2}\|\varphi_d\|_{L^2}.
\end{aligned}
\end{equation*}
Using \eqref{usefulest2}, we infer
\begin{equation*}
\begin{aligned}
&\left|\int_1^\infty\int_{\mathbb R^{d-1}}\nabla v^\sharp(y)\varphi_\diamond\varphi_d(y)dy'dy_d\right|\\
&\leq C\int_1^\infty\int_{\mathbb R^{d-1}}\left(\int_{\mathbb R^{d-1}\times\{0\}}\frac{1}{(y_d^2+|y'-\tilde{y}'|^2)^{d/2}}d\tilde{y}'\right)^{1/2}\\
&\qquad\qquad\qquad\left(\int_{\mathbb R^{d-1}\times\{0\}}\frac{1}{(y_d^2+|y'-\tilde{y}'|^2)^{d/2}}|v_0(\tilde{y}')|^2d\tilde{y}\right)^{1/2}|\varphi_\diamond\varphi_d(y)|dy'dy_d\\
&\leq C\|v_0\|_{L^2_{uloc}}\int_1^\infty\frac{1}{y_d}|\varphi_d(y_d)|dy_d\int_{\mathbb R^{d-1}}|\varphi_\diamond(y')|dy'\leq CR^{\frac{d-1}{2}}\|v_0\|_{L^2_{uloc}}\|\varphi_\diamond\|_{L^2}\|\varphi_d\|_{L^2}.
\end{aligned}
\end{equation*}

\subsection{Uniqueness}

By linearity, it is enough to prove uniqueness for $v=v(y)$ weak solution to
\begin{equation*}
\left\{
\begin{array}{rll}
-\nabla\cdot A(y)\nabla v&=0,&y_d>0,\\
v&=0,&y_d=0,
\end{array}
\right.
\end{equation*}
such that
\begin{equation}\label{estulocuniqflat}
\sup_{\xi\in\mathbb Z^{d-1}}\int_{\xi+(0,1)^{d-1}}\int_0^\infty|\nabla v|^2\leq C<\infty.
\end{equation}
Clearly, by Poincar\'e's inequality, for all $a>0$,
\begin{equation*}
\sup_{\xi\in\mathbb Z^{d-1}}\int_{\xi+(0,1)^{d-1}}\int_0^a|v|^2\leq Ca^2.
\end{equation*}
For $k\in\mathbb N$, we will take as a test function $\eta_k^2v\in H^1_0(\mathbb R^d_+)$ for an ad hoc cut-off $\eta_k$ such that $\eta_k\equiv 1$ on $(-k,k)^{d-1}\times (0,k)$ and $\supp\eta_k\subset(-k-1,k+1)^{d-1}\times(-1,k+1)$. 

We want to construct $\eta_k$ such that $\|\nabla\eta_k\|_{L^\infty}$ is bounded uniformly in $k$. Let $\eta\in C^\infty_c(B(0,1/2))$ such that $\int_{\mathbb R^d}\eta=1$. We define $\eta_k$ as folows
\begin{equation*}
\begin{aligned}
\eta_k(y)&:=\int_{\mathbb R^d}\mathbf{1}_{(-k-1/2,k+1/2)^{d-1}\times(-1/2,k+1/2)}(y-\tilde{y})\eta(\tilde{y})d\tilde{y}\\
&=\int_{(-k-1/2,k+1/2)^{d-1}\times(-1/2,k+1/2)}\eta(y-\tilde{y})d\tilde{y}.
\end{aligned}
\end{equation*}
For $y\in (-k,k)^{d-1}\times (0,k)$, $\supp(\eta(y-\cdot))\subset(-k-1/2,k+1/2)^{d-1}\times(-1/2,k+1/2)$, so $\eta_k(y)=1$. Moreover, $\supp\eta_k\subset(-k-1/2,k+1/2)^{d-1}\times(-1/2,k+1/2)+\supp\eta\subset(-k-1,k+1)^{d-1}\times(-1,k+1)$. Finally, convolution inequalities imply $\|\nabla\eta_k\|_{L^\infty}\leq\|\nabla\eta\|_{L^1}$.

Now, testing against $\eta_k^2v$, we get
\begin{equation*}
0=\int_{\mathbb R^d_+}A(y)\nabla v\cdot\nabla(\eta_k^2v)=\int_{\mathbb R^d_+}A(y)\eta_k^2\nabla v\cdot\nabla v+2\int_{\mathbb R^d_+}A(y)\eta_k\nabla v\cdot(\nabla\eta_k)v.
\end{equation*}
Therefore, letting 
\begin{equation*}
E_k:=\int_{(-k,k)^{d-1}\times(0,k)}|\nabla v|^2,
\end{equation*}
we have
\begin{equation}\label{estuniqflat}
E_k\leq C^*(E_{k+1}-E_k),
\end{equation}
where $C^*=C^*(d,N,\lambda,\|\nabla\eta\|_{L^1})$. Using the hole-filling trick, we get for fixed $k$ and for all $n\geq k$,
\begin{equation*}
E_k\leq \left(\frac{C^*}{C^*+1}\right)^{n-k}E_n.
\end{equation*}
Estimate \eqref{estulocuniqflat} implies $E_n\leq Cn^{d-1}$, so that
\begin{equation*}
E_k\leq C\left(\frac{C^*}{C^*+1}\right)^{n-k}n^{d-1}\stackrel{n\rightarrow\infty}{\longrightarrow}0,
\end{equation*}
and $E_k=0$. This concludes the  uniqueness proof.

\section{Estimates for a Dirichlet to Neumann operator}
\label{secDN}

The Dirichlet to Neumann operator $\DtoN$ is crucial in the proof of the well-posedness of the elliptic system in the bumpy half-space (see section \ref{secwpbumpyhp}). The key idea there is to carry out a domain decomposition. The Dirichlet to Neumann map is the tool enabling this domain decomposition. Since we are working in spaces of infinite energy to be useful $\DtoN$ has to be defined on $H^{1/2}_{uloc}$. Similar studies have been carried out in \cite{ABZ13_Kato} (context of water-waves), \cite{DGVNMnoslip} ($2$d Stokes system), \cite{DP_SC} ($3$d Stokes-Coriolis system).

We first define the Dirichlet to Neumann operator on $H^{1/2}(\mathbb R^{d-1})$:
\begin{equation*}
\DtoN:\ H^{1/2}(\mathbb R^{d-1})\longrightarrow \mathcal D',
\end{equation*}
such that for any $v_0\in H^{1/2}(\mathbb R^{d-1})$, for all $\varphi\in C^\infty_c(\mathbb R^{d-1})$,
\begin{equation*}
\langle\DtoN(v_0),\varphi\rangle_{\mathcal D',\mathcal D}:=\langle A(y)\nabla v\cdot e_d,\varphi\rangle_{\mathcal D',\mathcal D},
\end{equation*}
where $v$ is the unique weak solution to 
\begin{equation}\label{sysflatH1/2secdn}
\left\{\begin{array}{rll}
-\nabla\cdot A(y)\nabla v&=0,&y_d>0,\\
v&=v_0\in H^{1/2}(\mathbb R^{d-1}),&y_d=0.
\end{array}
\right.
\end{equation}

\begin{prop}\label{propdnh12}\  
\begin{enumerate}[label=(\arabic*)]
\item For all $\varphi\in C^\infty_c(\overline{\mathbb R^d_+})$, 
\begin{equation}\label{prop1dnh12}
\langle\DtoN(v_0),\varphi|_{y_d=0}\rangle_{\mathcal D',\mathcal D}=\langle A(y)\nabla v\cdot e_d,\varphi|_{y_d=0}\rangle_{\mathcal D',\mathcal D}=-\int_{\mathbb R^d_+}A(y)\nabla v\cdot\nabla\varphi.
\end{equation}
\item For all $\varphi\in C^\infty_c(\mathbb R^{d-1})$,
\begin{equation}\label{prop2dnh12}
\langle\DtoN(v_0),\varphi|_{y_d=0}\rangle_{\mathcal D',\mathcal D}=\int_{\mathbb R^{d-1}\times\{0\}}\int_{\mathbb R^{d-1}\times\{0\}}A(y)\nabla_yP(y,\tilde{y})\cdot e_dv_0(\tilde{y})d\tilde{y}\varphi(y)dy.
\end{equation}
\end{enumerate}
\end{prop}

For $y,\ \tilde{y}\in\mathbb R^{d-1}\times\{0\}$, let
\begin{equation*}
K(y,\tilde{y}):=A(y)\nabla_yP(y,\tilde{y})\cdot e_d
\end{equation*}
be the kernel appearing in \eqref{prop2dnh12}. Estimate \eqref{estnablaPbdary} of Proposition \ref{propestpeps} implies that
\begin{equation*}
|K(y,\tilde{y})|\leq \frac{C}{|y-\tilde{y}|^d},
\end{equation*}
for any $y,\ \tilde{y}\in \mathbb R^{d-1}\times\{0\}$, $y\neq\tilde{y}$ with $C=C(d,N,\lambda,[A]_{C^{0,\nu}})$.

Both formulas in Proposition \ref{propdnh12} follow from integration by parts. Because of \eqref{prop1dnh12}, it is clear that for all $v_0\in H^{1/2}(\mathbb R^{d-1})$, for all $\varphi\in C^{\infty}_c(\mathbb R^{d-1})$,
\begin{equation}\label{continuityestdnflathp}
\left|\langle\DtoN(v_0),\varphi\rangle\right|\leq C\|v_0\|_{H^{1/2}}\|\varphi\|_{H^{1/2}},
\end{equation}
with $C=C(d,N,\lambda)$, so that $\DtoN(v_0)$ extends as a continuous operator on $H^{1/2}(\mathbb R^{d-1})$. Another consequence of \eqref{prop1dnh12} is the following corollary.

\begin{cor}\label{cordnneg}
For all $v_0\in H^{1/2}(\mathbb R^{d-1})$, 
\begin{equation*}
\langle\DtoN(v_0),v_0\rangle=-\int_{\mathbb R^{d}_+}A(y)\nabla v\cdot\nabla v\leq 0,
\end{equation*}
where $v$ is the unique solution to \eqref{sysflatH1/2secdn}.
\end{cor}

Our next goal is to extend the definition of $\DtoN$ to $v_0\in H^{1/2}_{uloc}(\mathbb R^{d-1})$. We have to make sense of the duality product $\langle\DtoN(v_0),\varphi\rangle$. As for the definition of the solution to the flat half-space problem (see section \ref{secflathalfspace}), the basic idea is to use a cut-off function $\eta$ to split the definition between one part $\langle\DtoN(\eta v_0),\varphi\rangle$ where $\eta v_0\in H^{1/2}(\mathbb R^{d-1})$, and another part $\langle\DtoN((1-\eta)v_0),\varphi\rangle$ which does not see the singularity of the kernel $K(y,\tilde{y})$.

For $R>1$, there exists $\eta\in C^\infty_c(\mathbb R)$ such that
\begin{equation}\label{condeta}
0\leq\eta\leq 1,\quad\eta\equiv 1\ \mbox{on}\ (-R,R),\quad\supp\eta\subset(-R-1,R+1),\quad\|\eta'\|_{L^\infty}\leq 2.
\end{equation}
Let $v_0\in H^{1/2}_{uloc}(\mathbb R^{d-1})$. Let $R>1$ and $\varphi\in C^\infty_c(\mathbb R^{d-1})$ such that $\supp\varphi+B(0,1)\subset B(0,R)$. There exists $\eta\in C^\infty_c(\mathbb R)$ satisfying the conditions \eqref{condeta}. We define the action of $\DtoN(v_0)$ on $\varphi$ by
\begin{multline}\label{defdnh12ulocflat}
\langle\DtoN(v_0),\varphi\rangle_{\mathcal D',\mathcal D}:=\langle\DtoN(\eta(|\cdot|)v_0),\varphi\rangle_{H^{-1/2},H^{1/2}}\\
+\int_{\mathbb R^{d-1}\times\{0\}}\int_{\mathbb R^{d-1}\times\{0\}}K(y,\tilde{y})(1-\eta(|\tilde{y}'|))v_0(\tilde{y}')\varphi(y')d\tilde{y}dy.
\end{multline}
The fact that this definition does not depend on the cut-off $\eta\in C^\infty_c(\mathbb R)$ follows from Proposition \ref{propdnh12}.

The first term in the right-hand side of \eqref{defdnh12ulocflat} is estimated using \eqref{continuityestdnflathp} and the bound of Lemma \ref{lemh1/2H1/2uloc} between the $H^{1/2}$ norm of $\eta(|\cdot|)v_0$ and the $H^{1/2}_{uloc}$ norm of $v_0$. That yields
\begin{equation*}
\left|\langle\DtoN(\eta(|\cdot|)v_0),\varphi\rangle\right|\leq C\|\eta(|\cdot|)v_0\|_{H^{1/2}}\|\varphi\|_{H^{1/2}}\leq CR^{\frac{d-1}{2}}\|v_0\|_{H^{1/2}_{uloc}}\|\varphi\|_{H^{1/2}},
\end{equation*}
with $C=C(d,N,\lambda)$.

We deal with the integral part in the right hand side of \eqref{defdnh12ulocflat} in a way similar to the proof of estimates \eqref{usefulest1} and \eqref{usefulest2}. Using the fact that the supports of $(1-\eta(|y'|))v_0(y')$ on the one hand and $\varphi$ on the other hand are disjoint, we have
\begin{equation*}
\begin{aligned}
&\left|\int_{\mathbb R^{d-1}\times\{0\}}\int_{\mathbb R^{d-1}\times\{0\}}K(y,\tilde{y})(1-\eta(|\tilde{y}'|))v_0(\tilde{y}')\varphi(y')d\tilde{y}dy\right|\\
&\leq C\int_{\mathbb R^{d-1}\times\{0\}}\int_{\mathbb R^{d-1}\times\{0\}}\frac{1}{|y-\tilde{y}|^d}(1-\eta(|\tilde{y}'|))|v_0(\tilde{y}')||\varphi(y')|d\tilde{y}dy\\
&\leq C\int_{\mathbb R^{d-1}\times\{0\}}\left(\int_{\mathbb R^{d-1}\times\{0\}}\frac{1}{|y-\tilde{y}|^d}(1-\eta(|\tilde{y}'|))d\tilde{y}\right)^{1/2}\\
&\qquad\qquad\left(\int_{\mathbb R^{d-1}\times\{0\}}\frac{1}{|y-\tilde{y}|^d}(1-\eta(|\tilde{y}'|))|v_0(\tilde{y}')|^2d\tilde{y}\right)^{1/2}|\varphi(y')|dy\\
&\leq C\int_{\mathbb R^{d-1}\times\{0\}}\left(\int_1^\infty\frac{1}{r^2}dr\right)^{1/2}|\varphi(y')|dy\|v_0\|_{L^2_{uloc}}\\
&\leq CR^{\frac{d-1}{2}}\|v_0\|_{L^2_{uloc}}\|\varphi\|_{L^2},
\end{aligned}
\end{equation*}
with $C=C(d,N,\lambda,[A]_{C^{0,\nu}})$.

These results are put in a nutshell in the following proposition.

\begin{prop}\label{estdnflat}\ 
\begin{enumerate}[label=(\arabic*)]
\item For $v_0\in H^{1/2}(\mathbb R^{d-1})$, for any $\varphi\in C^\infty_c(\mathbb R^{d-1})$, we have
\begin{equation*}
\left|\langle\DtoN(v_0),\varphi\rangle\right|\leq C\|v_0\|_{H^{1/2}}\|\varphi\|_{H^{1/2}},
\end{equation*}
with $C=C(d,N,\lambda)$.
\item For $v_0\in H^{1/2}_{uloc}(\mathbb R^{d-1})$, for $R>1$ and any $\varphi\in C^\infty_c(\mathbb R^{d-1})$ such that 
\begin{equation*}
\supp\varphi+B(0,1)\subset B(0,R), 
\end{equation*}
we have
\begin{equation}\label{estdnflatest}
\left|\langle\DtoN(v_0),\varphi\rangle\right|\leq CR^{\frac{d-1}{2}}\|v_0\|_{H^{1/2}_{uloc}}\|\varphi\|_{H^{1/2}},
\end{equation}
with $C=C(d,N,\lambda,[A]_{C^{0,\nu}})$.
\end{enumerate}
\end{prop}

\section{Boundary layer corrector in a bumpy half-space}
\label{secwpbumpyhp}

This section is devoted to the well-posedness of the boundary layer problem
\begin{equation}\label{sysblbumpyhp}
\left\{
\begin{array}{rll}
-\nabla\cdot A(y)\nabla v&=0,&y_d>\psi(y'),\\
v&=v_0\in H^{1/2}_{uloc}(\mathbb R^{d-1}),&y_d=\psi(y'),
\end{array}
\right.
\end{equation}
in the bumpy half-space $\Omega_+:=\{y_d>\psi(y')\}$. For technical reasons, the boundary $\psi\in W^{1,\infty}(\mathbb R^{d-1})$ is assumed to be negative, i.e. $\psi(y')<0$ for all $y'\in\mathbb R^{d-1}$. We prove Theorem \ref{theoblbumpy} of the introduction which asserts the existence of a unique solution $v$ in the class
\begin{equation*}
\sup_{\xi\in\mathbb Z^{d-1}}\int_{\xi+(0,1)^{d-1}}\int_{\psi(y')}^\infty|\nabla v|^2dy_ddy'<\infty.
\end{equation*}

The idea is to split the bumpy half-space into two subdomains: a flat half-space $\mathbb R^d_+$ on the one hand and a bumpy channel $\Omega_\flat:=\{\psi(y')<y_d<0\}$ on the other hand. Both domains are connected by a transparent boundary condition involving the Dirichlet to Neumann operator $\DtoN$ defined in section \ref{secDN}. Therefore, solving \eqref{sysblhp} is equivalent to solving 
\begin{equation}\label{sysblbumpyc}
\left\{
\begin{array}{rll}
-\nabla\cdot A(y)\nabla v&=0,&0>y_d>\psi(y'),\\
v&=v_0\in H^{1/2}_{uloc}(\mathbb R^{d-1}),&y_d=\psi(y'),\\
A(y)\nabla v\cdot e_d&=\DtoN(v|_{y_d=0}),&y_d=0.
\end{array}
\right.
\end{equation}
This fact is stated in the following technical lemma.

\begin{lem}\label{lemeqsysdec}
If $v$ is a weak solution of \eqref{sysblbumpyc} in $\Omega_\flat$ such that
\begin{equation*}
\sup_{\xi\in\mathbb Z^{d-1}}\int_{\xi+(0,1)^{d-1}}\int_{\psi(y')}^0|\nabla v|^2dy_ddy'<\infty,
\end{equation*}
then $\tilde{v}$, defined by $\tilde{v}(y):=v(y)$ for $\psi(y')<y_d<0$ and $\tilde{v}|_{\mathbb R^d_+}$ is the unique solution to \eqref{sysblhp} with boundary condition $\tilde{v}|_{y_d=0^+}=v|_{y_d=0^-}$ given by Theorem \ref{theoblflat}, is a weak solution to \eqref{sysblbumpyhp}. Moreover, the reverse is also true. Namely, if $v$ is a weak solution to \eqref{sysblbumpyhp} in $\Omega_+$ such that 
\begin{equation*}
\sup_{\xi\in\mathbb Z^{d-1}}\int_{\xi+(0,1)^{d-1}}\int_{\psi(y')}^\infty|\nabla v|^2dy_ddy'<\infty,
\end{equation*}
then $v|_{\{\psi(y')<y_d<0\}}$ is a weak solution to \eqref{sysblbumpyc}.
\end{lem}

The main advantage of the domain decomposition is to make it possible to work in a channel, bounded in the vertical direction, in which one can rely on Poincar\'e type inequalities. Therefore our method is energy based, which makes it possible to deal with rough boundaries.

We now lift the boundary condition $v_0$. There exists $V_0$ such that 
\begin{equation*}
\sup_{\xi\in\mathbb Z^{d-1}}\int_{\xi+(0,1)^{d-1}}\int_{\psi(y')}^\infty|V_0|^2+|\nabla V_0|^2dy_ddy'\leq C\|v_0\|_{H^{1/2}_{uloc}}^2,
\end{equation*}
with $C=C(d,N,\|\psi\|_{W^{1,\infty}})$ and such that the trace of $V_0$ is $v_0$. Thus, $w:=v-V_0$ solves the system
\begin{equation}\label{sysblbumpychomo}
\left\{
\begin{array}{rll}
-\nabla\cdot A(y)\nabla w&=\nabla\cdot F,&0>y_d>\psi(y'),\\
w&=0,&y_d=\psi(y'),\\
A(y)\nabla w\cdot e_d&=\DtoN(w|_{y_d=0})+f,&y_d=0,
\end{array}
\right.
\end{equation}
where
\begin{equation*}
\begin{aligned}
F&:=A(y)\nabla V_0,\\
f&:=\DtoN(V_0|_{y_d=0})-A(y)\nabla V_0\cdot e_d.
\end{aligned}
\end{equation*}
Notice that the source terms satisfy the following estimates:
\begin{equation}\label{estFl2uloc}
\sup_{\xi\in\mathbb Z^{d-1}}\int_{\xi+(0,1)^{d-1}}\int_{\psi(y')}^0|F|^2dy_ddy'\leq C\|v_0\|_{H^{1/2}_{uloc}}^2,
\end{equation}
with $C=C(d,N,\lambda,\|\psi\|_{W^{1,\infty}})$ and for all $\varphi\in C^\infty_c(\mathbb R^{d-1})$ such that $B(0,R)\subset\supp\varphi\subset B(0,2R)$ for some $R>0$,
\begin{equation}\label{estfh1/2uloc}
\left|\langle f,\varphi\rangle\right|\leq CR^{\frac{d-1}{2}}\|v_0\|_{H^{1/2}_{uloc}}\|\varphi\|_{H^{1/2}},
\end{equation}
with $C=C(d,N,\lambda,[A]_{C^{0,\nu}},\|\psi\|_{W^{1,\infty}})$.

There are three steps in the proof of the well-posedness of \eqref{sysblbumpyhp}. Firstly, for $n\in\mathbb N$ we build approximate solutions $w_n=w_n(y)$ solving 
\begin{equation}\label{sysblbumpychomon}
\left\{
\begin{array}{rll}
-\nabla\cdot A(y)\nabla w_n&=\nabla\cdot F,&0>y_d>\psi(y'),\\
w_n&=0,&\{y_d=\psi(y')\}\cup\{|y'|=n\},\\
A(y)\nabla w_n\cdot e_d&=\DtoN(w_n|_{y_d=0})+f,&y_d=0,
\end{array}
\right.
\end{equation}
on $\Omega_{\flat,n}:=\{y'\in(-n,n)^{d-1},\ 0>y_d>\psi(y')\}$ and extend $w_n$ by $0$ on $\Omega_\flat\setminus\Omega_{\flat,n}$. We have that $w_n\in H^1(\Omega_\flat)$. This construction is utterly classical. Secondly, we aim at getting estimates uniform in $n$ on $w_n$ in the norm
\begin{equation}\label{sndstepbound}
\sup_{\xi\in\mathbb Z^{d-1}}\int_{\xi+(0,1)^{d-1}}\int_{\psi(y')}^0|\nabla w_n|^2dy_ddy'.
\end{equation}
This is done carrying out so-called Saint-Venant estimates in the bounded channel. We close this step by using a hole-filling argument. The method has been pioneered by Lady{\v{z}}enskaja and Solonnikov \cite{LS} for the Navier-Stokes system in a bounded channel. Here the situation is more involved because of the nonlocal operator $\DtoN$ on the upper boundary. The situation here is closer to \cite{DGVNMnoslip,DaGV11} ($2$d Stokes system) and \cite{DP_SC} ($3$d Stokes-Coriolis system). Finally, one has to check that weak limits of $w_n$ are indeed solutions of \eqref{sysblbumpychomo}. This step is straightforward because of the linearity of the equations. Uniqueness follows from the Saint-Venant estimate of the second step, with zero source terms.

We focus on the second step, which is by far the most intricate one. Let $r>0$, $y'_0\in\mathbb R^{d-1}$ and 
\begin{equation*}
\Omega_{\flat,y'_0,r}:=\{y'\in B(y'_0,r),\ 0>y_d>\psi(y')\}.
\end{equation*}
Let $w_r\in H^1(\Omega_\flat)$ be a weak solution to
\begin{equation}\label{sysblbumpychomor}
\left\{
\begin{array}{rll}
-\nabla\cdot A(y)\nabla w_r&=\nabla\cdot F_r,&0>y_d>\psi(y'),\\
w_r&=0,&y_d=\psi(y'),\\
A(y)\nabla w_r\cdot e_d&=\DtoN(w_r|_{y_d=0})+f_r,&y_d=0,
\end{array}
\right.
\end{equation}
such that $w_r=0$ on $\Omega_\flat\setminus\Omega_{\flat,y'_0,r}$, and where 
\begin{equation*}
F_r:=F\mathbf{1}_{\Omega_{\flat,y'_0,r}},\qquad f_r:=f\mathbf{1}_{B(y'_0,r)}.
\end{equation*}
Both $F_r$ and $f_r$ satisfy (respectively) the estimates \eqref{estFl2uloc} and \eqref{estfh1/2uloc} with constants uniform in $r$. Notice furthermore that $w_n$ defined above (see \eqref{sysblbumpychomon}) is equal to $w_r$ solution of \eqref{sysblbumpychomor} for $r:=n$ and $y'_0=0$.

For $k\in \mathbb N$, let 
\begin{equation*}
\Omega_{\flat,k}:=\{y'\in(-k,k)^{d-1},\ 0>y_d>\psi(y')\}.
\end{equation*}
Our goal is to estimate,
\begin{equation*}
E_k:=\int_{\Omega_{\flat,k}}|\nabla w_r|^2.
\end{equation*}
In the following, for $k,\ m\in\mathbb N$, $k,\ m\geq 1$, 
\begin{equation*}
\Sigma_k:=(-k,k)^{d-1},
\end{equation*}
and the set $\mathcal C_{k,m}$ denotes the family of cubes $T$ of volume $m^{d-1}$ contained in $\mathbb R^{d-1}\setminus\Sigma_{k+m-1}$ with vertices in $\mathbb Z^{d-1}$, i.e.
\begin{equation*}
\mathcal C_{k,m}:=\left\{T=\xi+(-m',m')^{d-1},\ \xi\in\mathbb Z^{d-1}\ \mbox{and}\ T\subset\mathbb R^{d-1}\setminus\Sigma_{k+m-1}\right\}.
\end{equation*}
Let also $\mathcal C_m$ be the family of all the cubes of volume $m^{d-1}$ with vertices in $\mathbb Z^{d-1}$
\begin{equation*}
\mathcal C_{m}:=\left\{T=\xi+(-m',m')^{d-1},\ \xi\in\mathbb Z^{d-1}\right\}.
\end{equation*}
Notice that for $k\geq \hat{k}\geq m'$,
\begin{equation*}
\mathcal C_{k,m}\subset\mathcal C_{\hat{k},m}\subset\mathcal C_{m',m}\subset\mathcal C_m.
\end{equation*}
For $T\in\mathcal C_{k,m}$, 
\begin{equation}\label{defomegaT}
E_T:=\int_{\Omega_{T}}|\nabla w_r|^2,\qquad \Omega_T:=\{y'\in T,\ 0>y_d>\psi(y')\}.
\end{equation}

\begin{prop}\label{propaprioriest}
There exists a constant $C^*=C^*(d,N,\lambda,[A]_{C^{0,\nu}},\|\psi\|_{W^{1,\infty}},\|v_0\|_{H^{1/2}_{uloc}})$ such that for all $r>0$, $y'_0\in\mathbb R^{d-1}$, for all $k,\ m\in\mathbb N$, $m\geq 3$ and $k\geq m/2=m'$, for any weak solution $w_r\in H^1(\Omega_\flat)$ of \eqref{sysblbumpychomor}, the following bound holds
\begin{equation}\label{stvenantest}
E_k\leq C^*\left(k^{d-1}+E_{k+m}-E_k+\frac{k^{3d-5}}{m^{3d-3}}\sup_{T\in\mathcal C_{k,m}}E_T\right).
\end{equation}
Notice that $C^*$ is independent of $r$ and $y'_0$.
\end{prop}

The crucial point for the control of the large-scale energies in \eqref{stvenantest} is the fact that the power $3d-5$ of $k$ is strictly smaller that the power $3d-3$ of $m$. Before tackling the proof of Proposition \ref{propaprioriest}, let us explain how to infer from \eqref{stvenantest} an a priori bound uniform in $n$ on $w_n$ solution of \eqref{sysblbumpychomon}.

\subsection{Proof of the a priori bound}
\label{secproofaprioribound}

Let $C^*$ be given by Proposition \ref{propaprioriest}, and let
\begin{equation}\label{eqAB}
A:=\sum_{k=1}^\infty\left(\frac{C^*}{C^*+1}\right)^k(2k-1)^{d-1}<\infty,\qquad B:=\sum_{k=1}^\infty\left(\frac{C^*}{C^*+1}\right)^k(2k-1)^{3d-5}<\infty.
\end{equation}
We now choose an integer $m$ so that 
\begin{equation}\label{eqdefm}
m\geq 3,\quad m\ \mbox{is even}\quad\mbox{and}\quad 1-2^{5-3d}\frac{B}{m^2}>\frac{1}{2}.
\end{equation}
Notice that $m=m(d,N,\lambda,[A]_{C^{0,\nu}},\|\psi\|_{W^{1,\infty}},\|v_0\|_{H^{1/2}_{uloc}})$, but is independent of $r$ and $y_0'$. The reason for taking $m$ even is technical; it is only used in the translation argument below.

Take $n=lm=2lm'$, with $l\in\mathbb N$, $l\geq 1$, and take $w_n$ to be the solution of \eqref{sysblbumpychomon}. There exists $T^*\in \mathcal C_{m}$ such that $T^*\subset\Sigma_n$ and $E_{T^*}=\sup_{T\in\mathcal C_{m}}E_T$. By definition, there is $\xi^*\in\mathbb Z^{d-1}$ for which $T^*=\xi^*+(-m',m')^{d-1}$. We want to center $T^*$ at zero by simply translating the origin. Doing so, $w^*_n(y):=w_n(y'+\xi^*,y_d)$ is a solution of \eqref{sysblbumpychomor} with $y'_0:=-\xi^*,\ r=n$ and
\begin{equation*}
\begin{aligned}
A^*(y):=A(y'+\xi^*,y_d),&\quad \psi^*(y'):=\psi(y'+\xi^*),\quad v_0^*(y'):=v_0(y'+\xi^*),\\
F^*(y):=F(y'+y^*,y_d)&\quad\mbox{and}\quad f^*(y'):=f(y'+\xi^*).
\end{aligned}
\end{equation*}
Notice that 
\begin{equation*}
[A^*]_{C^{0,\mu}}=[A]_{C^{0,\nu}},\quad \|\psi^*\|_{W^{1,\infty}}=\|\psi\|_{W^{1,\infty}}\quad\mbox{and}\quad \|v^*_0\|_{H^{1/2}_{uloc}}=\|v_0\|_{H^{1/2}_{uloc}},
\end{equation*}
so that $w^*_n$ satisfies the Saint-Venant estimate \eqref{stvenantest} with the same constant $C^*$. Furthermore, $E_{m'}=E_{T^*}$.

\begin{lem}\label{lemdownwardinduction}
We have the following a priori bound
\begin{equation*}
E_{m'}\leq 2^{2-d}Am^{d-1},
\end{equation*}
where $A$ is defined by \eqref{eqAB}.
\end{lem}

The Lemma is obtained by downward induction, using a hole-filling type argument. Since $w^*_n$ is supported in $\Omega_{\flat,2n}$, we start from $k$ sufficiently large in \eqref{stvenantest}. For $k=2n+m'=(4l+1)m'$, estimate \eqref{stvenantest} implies
\begin{equation*}
E_{(4l+1)m'}\leq \frac{C^*}{C^*+1}((4l+1)m')^{d-1},
\end{equation*}
because $E_T=0$ for any $T\in\mathcal C_{(4l+1)m',m}$. Then,
\begin{equation*}
E_{(2(2l-1)+1)m'}=E_{(4l+1)m'-m}\leq \frac{C^*}{C^*+1}(2(2l-1)+1)^{d-1}(m')^{d-1}+\left(\frac{C^*}{C^*+1}\right)^2(4l+1)^{d-1}(m')^{d-1}.
\end{equation*}
Let $p\in\{0,\ldots\ 2l-1\}$. We then have
\begin{multline*}
E_{(2p+1)m'}\leq \frac{C^*}{C^*+1}(2p+1)^{d-1}(m')^{d-1}+\ldots\ \left(\frac{C^*}{C^*+1}\right)^{2l-p}(4l+1)^{d-1}(m')^{d-1}\\
+\frac{2^{5-3d}}{m^2}\left[\frac{C^*}{C^*+1}(2p+1)^{3d-5}+\ldots\ \left(\frac{C^*}{C^*+1}\right)^{2l-p}(4l+1)^{3d-5}\right]E_{m'}.
\end{multline*}
Eventually, for $p=0$ 
\begin{multline*}
E_{m'}\leq \frac{C^*}{C^*+1}(m')^{d-1}+\left(\frac{C^*}{C^*+1}\right)^2(3m')^{d-1}+\ldots\ \left(\frac{C^*}{C^*+1}\right)^{2l-p}(4l+1)^{d-1}(m')^{d-1}\\
+\frac{2^{5-3d}}{m^2}\left[\frac{C^*}{C^*+1}+\ldots\ \left(\frac{C^*}{C^*+1}\right)^{2l-p}(4l+1)^{3d-5}\right]E_{m'}\leq 2^{1-d}Am^{d-1}+\frac{2^{5-3d}}{m^2}BE_{m'}.
\end{multline*}
Therefore, 
\begin{equation*}
\frac{E_{m'}}{2}<\left(1-\frac{B}{m^2}\right)E_{m'}\leq 2^{1-d}Am^{d-1},
\end{equation*}
which proves Lemma \ref{lemdownwardinduction}.

Finally, 
\begin{equation*}
\sup_{\xi\in\mathbb Z^{d-1}}\int_{\xi+(0,1)^{d-1}}\int_{\psi(y')}^0|\nabla w_n|^2\leq E_m\leq 2^{2-d}Am^{d-1},
\end{equation*}
which proves the a priori bound in the norm \eqref{sndstepbound} uniformly in $n$. 

\subsection{Proof of Proposition \ref{propaprioriest}}
\label{secproofpropaprioriest}

\subsubsection*{Construction of a cut-off}

Let $\eta\in C^\infty(B(0,1/2))$ such that $\eta\geq 0$ and $\int_{\mathbb R^d}\eta=1$. For all $k\in\mathbb N$, let $\eta_k=\eta_k(y')$ be defined by
\begin{equation*}
\eta_k(y')=\int_{\mathbb R^{d-1}}\mathbf{1}_{[-k-1/2,k+1/2]^{d-1}}(y'-\tilde{y}')\eta(\tilde{y}')d\tilde{y}'=\int_{[-k-1/2,k+1/2]^{d-1}}\eta(y'-\tilde{y}')d\tilde{y}'.
\end{equation*}
For all $k\in\mathbb N$, we have the following properties:
\begin{equation*}
\eta_k\equiv 1\ \mbox{on}\ [-k,k]^{d-1},\quad \supp\eta_k\subset[-k-1,k+1]^{d-1},\quad \eta_k\in C^\infty_c(\mathbb R^{d-1})
\end{equation*}
and most importantly, we have the control
\begin{equation*}
\|\nabla\eta_k\|_{L^\infty}\leq\|\nabla\eta\|_{L^1}
\end{equation*}
uniformy in $k$.

\subsubsection*{Energy estimate}

Testing the system \eqref{sysblbumpychomor} against $\eta_k^2w_r$ we get
\begin{multline}\label{eqtestetak2}
\int_{\Omega_\flat}\eta_k^2A(y)\nabla w_r\cdot\nabla w_r=-\int_{\Omega_\flat}2\eta_kA(y)\nabla w_r\cdot\nabla\eta_k w_r\\+\langle\nabla\cdot F,\eta_k^2w_r\rangle+\langle f,\eta_k^2w_r\rangle+\langle\DtoN(w_r|_{y_d=0}),\eta_k^2w_r\rangle.
\end{multline}
By ellipticity, we have
\begin{equation*}
\lambda\int_{\Omega_\flat}\eta_k^2|\nabla w_r|^2\leq\int_{\Omega_\flat}\eta_k^2A(y)\nabla w_r\cdot\nabla w_r.
\end{equation*}
The following estimate (or variations of it) is of constant use: by the trace theorem and Poincar\'e inequality
\begin{equation}\label{eststvenanta}
\begin{aligned}
\left(\int_{\Sigma_{k+1}}\eta_k^4|w_r(y',0)|dy'\right)^{1/2}&\leq \|\eta_k^2w_r\|_{H^{1/2}}\leq C\|\eta_k^2w_r\|_{H^1(\Omega_\flat)}\leq C\|\nabla(\eta_k^2w_r)\|_{L^2(\Omega_\flat)}\\
&\leq C(E_{k+1}-E_k)^{1/2}+C'\left(\int_{\Omega_\flat}\eta_k^4|\nabla w|^2\right)^{1/2},
\end{aligned}
\end{equation}
with $C=C(d,\|\psi\|_{W^{1,\infty}},\|\eta\|_{L^1})$ and $C'=C'(d)$. We now estimate every term on the right hand side of \eqref{eqtestetak2}. We have,
\begin{equation*}
\begin{aligned}
\left|\int_{\Omega_\flat}2\eta_kA(y)\nabla w_r\cdot\nabla\eta_k w_r\right|&\leq\frac{2}{\lambda}\left(\int_{\Omega_\flat}\eta_k^2|\nabla w_r|^2\right)^{1/2}\left(\int_{\Omega_\flat}|\nabla\eta_k|^2|w_r|^2\right)^{1/2}\\
&\leq C\left(\int_{\Omega_\flat}\eta_k^2|\nabla w_r|^2\right)^{1/2}(E_{k+1}-E_k)^{1/2},
\end{aligned}
\end{equation*}
with $C=C(\lambda,\|\eta\|_{L^1})$. We also have,
\begin{equation*}
\begin{aligned}
|\langle\nabla\cdot F,\eta^2_kw_r\rangle|&=|\langle F,\nabla(\eta_k^2w_r)\rangle|\\
&\leq Ck^{\frac{d-1}{2}}(E_{k+1}-E_k)^{1/2}+C'k^{\frac{d-1}{2}}\left(\int_{\Omega_\flat}\eta_k^4|\nabla w_r|^2\right)^{1/2},
\end{aligned}
\end{equation*}
where $C=C(\|v_0\|_{H^{1/2}_{uloc}},\|\nabla\eta\|_{L^1})$ and $C'=C'(\|v_0\|_{H^{1/2}_{uloc}})$, and by the trace theorem and Poincar\'e inequality
\begin{equation*}
\begin{aligned}
|\langle f,\eta_k^2w_r\rangle|&\leq Ck^\frac{d-1}{2}\|\eta_k^2w_r\|_{H^{1/2}}\leq Ck^\frac{d-1}{2}\|\nabla(\eta_k^2w_r)\|_{L^2}\\
&\leq Ck^\frac{d-1}{2}(E_{k+1}-E_k)^{1/2}+C'k^{\frac{d-1}{2}}\left(\int_{\Omega_\flat}\eta_k^4|\nabla w_r|^2\right)^{1/2},
\end{aligned}
\end{equation*}
with $C=C(d,\|\psi\|_{W^{1,\infty}},\|v_0\|_{H^{1/2}_{uloc}},\|\nabla\eta\|_{L^1})$ and $C'=C'(\|v_0\|_{H^{1/2}_{uloc}})$. We have now to tackle the non local term involving the Dirichlet to Neumann operator. We split this term into
\begin{multline*}
\langle\DtoN(w_r|_{y_d=0}),\eta_k^2w_r\rangle=\langle\DtoN((1-\eta_{k+m-1}^2)w_r|_{y_d=0}),\eta_k^2w_r\rangle\\
+\langle\DtoN((\eta_{k+m-1}^2-\eta_k^2)w_r|_{y_d=0}),\eta_k^2w_r\rangle+\langle\DtoN(\eta_{k}^2w_r),\eta_k^2w_r\rangle.
\end{multline*}
By Corollary \ref{cordnneg}, 
\begin{equation*}
\langle\DtoN(\eta_{k}^2w_r),\eta_k^2w_r\rangle\leq 0.
\end{equation*}
Relying on Proposition \ref{estdnflat} and on estimate \eqref{estdnflatest}, we get
\begin{equation*}
\begin{aligned}
&|\langle\DtoN((\eta_{k+m-1}^2-\eta_k^2)w_r|_{y_d=0}),\eta_k^2w_r\rangle|\leq Ck^{\frac{d-1}{2}}\|(\eta_{k+m-1}^2-\eta_k^2)w_r|_{y_d=0}\|_{H^{1/2}_{uloc}}\|\eta_k^2w_r\|_{H^{1/2}}\\
&\leq C(E_{k+m}-E_k)^{1/2}\left(\int_{\Omega_\flat}\eta_k^4|\nabla w_r|^2\right)^{1/2}+C(E_{k+m}-E_k)^{1/2}(E_{k+1}-E_k)^{1/2},
\end{aligned}
\end{equation*}
with $C=C(d,N,\lambda,[A]_{C^{0,\nu}},\|\psi\|_{W^{1,\infty}},\|\eta\|_{L^1})$. Notice that the bound \eqref{continuityestdnflathp} for the Dirichlet to Neumann operator in $H^{1/2}$ here is actually enough, since $w_r$ is compactly supported. However, when dealing with solutions not compactly supported, as for the uniqueness proof in section \ref{secendproofbumpy}, we have to use the result of Proposition \ref{estdnflat}.

\subsubsection*{Control of the non local term}

\begin{lem}\label{lemcontrolnonlocalterm}
For all $m\geq 3$, all $k\geq m'=m/2$, we have
\begin{equation}\label{lemcontrolnonlocaltermest}
\int_{\Sigma_{k+1}}\left(\int_{\mathbb R^{d-1}}\frac{1}{|y'-\tilde{y}'|^d}(1-\eta^2_{k+m-1})|w_r(\tilde{y}',0)|d\tilde{y}'\right)^2dy'\leq C\frac{k^{3d-5}}{m^{3d-3}}\sup_{T\in\mathcal C_{k,m}}E_T,
\end{equation}
where $C=C(d)$.
\end{lem}

Let $y'\in \Sigma_{k+1}$ be fixed. We have
\begin{equation*}
\begin{aligned}
&\int_{\mathbb R^{d-1}}\frac{1}{|y'-\tilde{y}'|^d}(1-\eta^2_{k+m-1})|w_r(\tilde{y}',0)|d\tilde{y}'\\
&=\sum_{j=1}^\infty\int_{\mathbb R^{d-1}}\frac{1}{|y'-\tilde{y}'|^d}(\eta^2_{k+(j+1)(m-1)}-\eta^2_{k+m-1})|w_r(\tilde{y}',0)|d\tilde{y}'\\
&=\sum_{j=1}^\infty\int_{\Sigma_{k+(j+1)(m-1)+1}\setminus\Sigma_{k+j(m-1)}}\frac{1}{|y'-\tilde{y}'|^d}|w_r(\tilde{y}',0)|d\tilde{y}'\\
&=\sum_{j=1}^\infty\sum_{T\in\mathcal C_{k,j,m}}\int_{T}\frac{1}{|y'-\tilde{y}'|^d}|w_r(\tilde{y}',0)|d\tilde{y}',
\end{aligned}
\end{equation*}
where $\mathcal C_{k,j,m}$ is a family of disjoint cubes $T=\xi+(-m',m')^{d-1}$ such that $T\subset\Sigma_{k+(j+1)(m-1)+1}\setminus\Sigma_{k+j(m-1)}$ and
\begin{equation*}
\bigsqcup_{T\in \mathcal C_{k,j,m}}T=\Sigma_{k+(j+1)(m-1)+1}\setminus\Sigma_{k+j(m-1)}.
\end{equation*}
For all $T\in \mathcal C_{k,j,m}$, by Cauchy-Schwarz, trace theorem and Poincar\'e inequality
\begin{equation*}
\begin{aligned}
\int_{T}\frac{1}{|y'-\tilde{y}'|^d}|w_r(\tilde{y}',0)|d\tilde{y}'&\leq\left(\int_T\frac{1}{|y'-\tilde{y}'|^{2d}}d\tilde{y}'\right)^{1/2}\left(\int_T|w(\tilde{y}',0)|^2d\tilde{y}'\right)^{1/2}\\
&\leq C\left(\int_T\frac{1}{|y'-\tilde{y}'|^{2d}}d\tilde{y}'\right)^{1/2}\left(\int_{\Omega_T}|\nabla w|^2d\tilde{y}'\right)^{1/2}\\
&\leq C\left(\int_T\frac{1}{|y'-\tilde{y}'|^{2d}}d\tilde{y}'\right)^{1/2}\left(\sup_{T\in \mathcal C_{k,j,m}}E_T\right)^{1/2},
\end{aligned}
\end{equation*}
where $\Omega_T$ and $E_T$ are defined in \eqref{defomegaT}. Notice that the constant $C$ in the last inequality only depends on $d$ and on $\|\psi\|_{W^{1,\infty}}$. Moreover, for any $T\in \mathcal C_{k,j,m}$, 
\begin{equation*}
\left(\int_T\frac{1}{|y'-\tilde{y}'|^{2d}}d\tilde{y}'\right)^{1/2}\leq \frac{m^{\frac{d-1}{2}}}{(k+j(m-1)-|y'|)^d},
\end{equation*}
and the number of elements of $\mathcal C_{k,j,m}$ is bounded by
\begin{equation*}
\#\mathcal C_{k,j,m}=\frac{\left|\Sigma_{k+(j+1)(m-1)+1}\setminus\Sigma_{k+j(m-1)}\right|}{m^{d-1}}\lesssim\frac{(k+j(m-1))^{d-2}}{m^{d-2}}. 
\end{equation*}
Therefore,
\begin{equation*}
\begin{aligned}
&\int_{\mathbb R^{d-1}}\frac{1}{|y'-\tilde{y}'|^d}(1-\eta^2_{k+m-1})|w_r(\tilde{y}',0)|d\tilde{y}'\\
&\leq C\left(\sup_{T\in \mathcal C_{k,j,m}}E_T\right)^{1/2}\sum_{j=1}^\infty\sum_{T\in\mathcal C_{k,j,m}}\frac{m^{\frac{d-1}{2}}}{(k+j(m-1)-|y'|)^d}\\
&\leq C\left(\sup_{T\in \mathcal C_{k,j,m}}E_T\right)^{1/2}\sum_{j=1}^\infty\frac{1}{m^{\frac{d-3}{2}}}\frac{(k+j(m-1))^{d-2}}{(k+j(m-1)-|y'|)^d}\\
&\leq C\left(\sup_{T\in \mathcal C_{k,j,m}}E_T\right)^{1/2}\frac{(k+m-1)^{d-2}}{m^{\frac{d-1}{2}}(k+m-1-|y'|)^{d-1}},
\end{aligned}
\end{equation*}
with $C=C(d)$. Eventually, we get for $m\geq 3$
\begin{equation*}
\begin{aligned}
&\int_{\Sigma_{k+1}}\left(\int_{\mathbb R^{d-1}}\frac{1}{|y'-\tilde{y}'|^d}(1-\eta^2_{k+m-1})|w_r(\tilde{y}',0)|d\tilde{y}'\right)^2dy'\\
&\leq C\left(\sup_{T\in \mathcal C_{k,j,m}}E_T\right)\frac{(k+m-1)^{2d-4}}{m^{d-1}}\int_{\Sigma_{k+1}}\frac{1}{(k+m-1-|y'|)^{2d-2}}dy'\\
&\leq C\left(\sup_{T\in \mathcal C_{k,j,m}}E_T\right)\frac{(k+m-1)^{2d-4}}{m^{d-1}}\frac{(k+1)^{d-1}}{(m-2)^{2d-2}}\leq C\frac{k^{3d-5}}{m^{3d-3}}\sup_{T\in \mathcal C_{k,j,m}}E_T,
\end{aligned}
\end{equation*}
with $C=C(d)$, the last inequality being only true on condition that $k\geq m/2=m'$. This proves Lemma \ref{lemcontrolnonlocalterm}.

In particular, by the definition of $\DtoN$ in \eqref{defdnh12ulocflat}, by the fact that $(1-\eta^2_{k+m-1})w_r(\tilde{y}',0)$ and $\eta_k^2w_r(y',0)$ have disjoint support, by estimate \eqref{lemcontrolnonlocaltermest} and by the bound \eqref{eststvenanta} we get
\begin{equation*}
\begin{aligned}
&|\langle\DtoN((1-\eta_{k+m-1}^2)w_r|_{y_d=0}),\eta_k^2w_r\rangle|\\
&\leq C\left|\int_{\mathbb R^{d-1}}\int_{\mathbb R^{d-1}}\frac{1}{|y'-\tilde{y}'|^d}(1-\eta^2_{k+m-1}(\tilde{y}'))|w_r(\tilde{y}',0)|\eta_k^2(y')|w_r(y',0)|d\tilde{y}'dy'\right|\\
&\leq C\left(\int_{\Sigma_{k+1}}\eta_k^4|w_r(y',0)|^2dy'\right)^{1/2}\left(\int_{\Sigma_{k+1}}\left(\int_{\mathbb R^{d-1}}\frac{1-\eta^2_{k+m-1}(\tilde{y}')}{|y'-\tilde{y}'|^d}|w_r(y',0)|d\tilde{y}'\right)^2dy'\right)^{1/2},\\
&\leq C\frac{k^\frac{3d-5}{2}}{m^\frac{3d-3}{2}}\left(\int_{\Sigma_{k+1}}\eta_k^4|w_r(y',0)|^2dy'\right)^{1/2}\left(\sup_{T\in \mathcal C_{k,j,m}}E_T\right)^\frac{1}{2}\\
&\leq C\frac{k^\frac{3d-5}{2}}{m^\frac{3d-3}{2}}\left[(E_{k+1}-E_k)^{1/2}+\left(\int_{\Omega_\flat}\eta_k^4|\nabla w|^2\right)^{1/2}\right]\left(\sup_{T\in \mathcal C_{k,j,m}}E_T\right)^\frac{1}{2},
\end{aligned}
\end{equation*}
with $C=C(d,N,\lambda,[A]_{C^{0,\nu}},\|\psi\|_{W^{1,\infty}})$.

\subsubsection*{End of the proof of the Saint-Venant estimate}

Combining all our bounds and using 
\begin{equation*}
E_{k+1}-E_k\leq E_{k+m}-E_k,\qquad \eta_k^4\leq C(\|\eta\|_{L^\infty})\eta_k^2
\end{equation*}
whenever possible, we get from \eqref{eqtestetak2} the following estimate
\begin{multline*}
\lambda\int_{\Omega_\flat}\eta_k^2|\nabla w_r|^2\leq C\left(\int_{\Omega_\flat}\eta_k^2|\nabla w_r|^2\right)^{1/2}(E_{k+m}-E_k)^{1/2}+Ck^\frac{d-1}{2}(E_{k+m}-E_k)^{1/2}+C(E_{k+m}-E_k)\\
+Ck^\frac{d-1}{2}\left(\int_{\Omega_\flat}\eta_k^2|\nabla w_r|^2\right)^{1/2}+C(E_{k+m}-E_k)^{1/2}\left(\int_{\Omega_\flat}\eta_k^2|\nabla w_r|^2\right)^{1/2}\\
+C\frac{k^\frac{3d-5}{2}}{m^\frac{3d-3}{2}}\left[(E_{k+1}-E_k)^{1/2}+\left(\int_{\Omega_\flat}\eta_k^2|\nabla w|^2\right)^{1/2}\right]\left(\sup_{T\in \mathcal C_{k,j,m}}E_T\right)^{1/2},
\end{multline*}
with $C=C(d,N,\lambda,[A]_{C^{0,\nu}},\|v_0\|_{H^{1/2}_{uloc}},\|\psi\|_{W^{1,\infty}})$. Swallowing every term of the type
\begin{equation*}
\int_{\Omega_\flat}\eta_k^2|\nabla w|^2
\end{equation*}
in the left hand side, we end up with the Saint-Venant estimate \eqref{stvenantest}. This concludes the proof of Proposition \ref{propaprioriest}.

\subsection{End of the proof of Theorem \ref{theoblbumpy}}
\label{secendproofbumpy}

Extracting subsequences using a classical diagonal argument and passing to the limit in the weak formulation of \eqref{sysblbumpychomon} relying on the continuity of the Dirichlet to Neumann map asserted in estimate \eqref{continuityestdnflathp} yields the existence of a weak solution $w$ to the system \eqref{sysblbumpychomo}. In addition, the weak solution satisfies the bound 
\begin{equation}\label{estwuloc}
\sup_{\xi\in\mathbb Z^{d-1}}\int_{\xi+(0,1)^{d-1}}\int_{\psi(y')}^\infty|\nabla w|^2dy_ddy'\leq 2^{2-d}Am^{d-1}<\infty.
\end{equation}

Let us turn to the uniqueness of the solution to \eqref{sysblbumpychomo} satisfying the bound \eqref{estwuloc}. By linearity of the problem, it is enough to prove the uniqueness for zero source terms. Assume $w\in H^1_{loc}(\Omega_\flat)$ is a weak solution to \eqref{sysblbumpychomo} with $f=F=0$ satisfying
\begin{equation}\label{boundulocC_0}
\sup_{\xi\in\mathbb Z^{d-1}}\int_{\xi+(0,1)^{d-1}}\int_{\psi(y')}^0|\nabla w|^2\leq C_0<\infty.
\end{equation}
Repeating the estimates leading to Proposition \ref{propaprioriest} (see section \ref{secproofpropaprioriest}), we infer that for the same constant $C^*$ appearing in the Saint-Venant estimate \eqref{stvenantest} and for $m$ defined by \eqref{eqdefm}, for $k\in\mathbb N$, $k\geq m/2=m'$,
\begin{equation}\label{stvenantestzerosource}
E_k\leq C^*\left(E_{k+m}-E_k+\frac{k^{3d-5}}{m^{3d-3}}\sup_{T\in\mathcal C_{k,m}}E_T\right).
\end{equation}
The fact that $w$, unlike $w_n$, does not vanish outside $\Omega_{\flat,n}$ does not lead to any difference in the proof of this estimate.

Since 
\begin{equation*}
\sup_{T\in\mathcal C_{m}}E_T<\infty,
\end{equation*}
for any $\varepsilon$, there exists $T^*_\varepsilon\in\mathcal C_{m}$ such that 
\begin{equation}\label{eqsqueeze}
\sup_{T\in\mathcal C_{m}}E_T-\varepsilon\leq E_{T^*_\varepsilon}\leq \sup_{T\in\mathcal C_{m}}E_T.
\end{equation}
Again, $T^*_\varepsilon:=\xi^*_\varepsilon+(-m',m')^{d-1}$ for $\xi^*_\varepsilon\in \mathbb Z^{d-1}$, and we can translate $T^*_\varepsilon$ so that it is centered at the origin as has been done in section \ref{secproofaprioribound}. Estimate \eqref{stvenantestzerosource} still holds. For any $n\in\mathbb N$, $E_n\leq C_0n^{d-1}$ where $C_0$ is defined by \eqref{boundulocC_0}. The idea is now to carry out a downward iteration. For any $n=(2l+1)m'$ with $l\in\mathbb N$, $l\geq 1$ fixed, for $p\in\{1,\ldots\ l-1\}$ one can show that
\begin{equation*}
\begin{aligned}
E_{(2p+1)m'}&\leq \left[\frac{C^*}{C^*+1}+\left(\frac{C^*}{C^*+1}\right)^2+\ldots\ \left(\frac{C^*}{C^*+1}\right)^{l-p}\right]E_n\\
&\qquad\qquad+\frac{2^{5-3d}}{m^2}\left[\frac{C^*}{C^*+1}(2p+1)^{3d-5}+\ldots\ \left(\frac{C^*}{C^*+1}\right)^{l-p}(2l+1)^{3d-5}\right]\sup_{T\in\mathcal C_{m}}E_T\\
&\leq C_0\frac{C^*+1}{2C^*+1}\left(\frac{C^*}{C^*+1}\right)^{l-p+1}n^{d-1}\\
&\qquad\qquad+\frac{2^{5-3d}}{m^2}\left[\frac{C^*}{C^*+1}(2p+1)^{3d-5}+\ldots\ \left(\frac{C^*}{C^*+1}\right)^{l-p}(2l+1)^{3d-5}\right]\sup_{T\in\mathcal C_{m}}E_T.
\end{aligned}
\end{equation*}
Thus,
\begin{equation*}
\begin{aligned}
E_{m'}&\leq C_0\frac{C^*+1}{2C^*+1}\left(\frac{C^*}{C^*+1}\right)^{l}(2l+1)^{d-1}(m')^{d-1}+\frac{2^{5-3d}}{m^2}B\sup_{T\in\mathcal C_{m}}E_T\\
&\leq C_0\frac{C^*+1}{2C^*+1}\left(\frac{C^*}{C^*+1}\right)^{2l+1}(2l+1)^{d-1}(m')^{d-1}+\frac{2^{5-3d}}{m^2}B(E_{m'}+\varepsilon).
\end{aligned}
\end{equation*}
From this we infer using \eqref{eqdefm} that
\begin{equation*}
E_{m'}\leq 2C_0\frac{C^*+1}{2C^*+1}\left(\frac{C^*}{C^*+1}\right)^{2l+1}(2l+1)^{d-1}(m')^{d-1}+\frac{2^{6-3d}}{m^2}B\varepsilon\stackrel{l\rightarrow\infty}{\longrightarrow}\frac{2^{6-3d}}{m^2}B\varepsilon.
\end{equation*}
Therefore, from equation \eqref{eqsqueeze}
\begin{equation*}
\sup_{T\in\mathcal C_{m}}E_T\leq \left(1+\frac{2^{6-3d}}{m^2}B\right)\varepsilon,
\end{equation*}
which eventually leads to $\sup_{T\in\mathcal C_{m}}E_T=0$, or in other words $w=0$.

Combining this existence and uniqueness result for the system \eqref{sysblbumpychomo} in the bumpy channel $\Omega_\flat$ with Lemma \ref{lemeqsysdec} and Theorem \ref{theoblflat} about the well-posedness in the flat half-space finishes the proof of Theorem \ref{theoblbumpy}.

\section{Improved regularity over Lipschitz boundaries}
\label{secimpro}

The goal in this section is to prove Theorem \ref{theolipdowmicro} of the introduction. Let us recall the result we prove in the following proposition.

\begin{prop}\label{proplipdowmicro}
For all $\nu>0$, $\gamma>0$, there exists $C>0$ and $\varepsilon_0>0$ such that for all $\psi\in W^{1,\infty}(\mathbb R^{d-1})$, $-1<\psi<0$ and $\|\nabla\psi\|_{L^\infty}\leq\gamma$, for all $A\in\mathcal A^\nu$, for all $0<\varepsilon<(1/2)\varepsilon_0$, for all weak solution $u^\varepsilon$ to \eqref{sysuepsintro}, for all $r\in[\varepsilon/\varepsilon_0,1/2]$
\begin{equation}\label{estlipdownmicroinprop}
\intbar_{D^\varepsilon_\psi(0,r)}|\nabla u^\varepsilon|^2\leq C\intbar_{D^\varepsilon_\psi(0,1)}|\nabla u^\varepsilon|^2,
\end{equation}
or equivalently,
\begin{equation*}
\intbar_{D^\varepsilon_\psi(0,r)}|u^\varepsilon|^2\leq Cr^2\intbar_{D^\varepsilon_\psi(0,1)}|u^\varepsilon|^2,
\end{equation*}
with $C=C(d,N,\lambda,\nu,\gamma,[A]_{C^{0,\nu}})$.
\end{prop}

We rely on a compactness argument inspired by the pioneering work of Avellaneda and Lin \cite{alin,Alin90P}, and our recent work \cite{BLtailosc}. The proof is in two steps. Firstly, we carry out the compactness argument. Secondly, we iterate the estimate obtained in the first step, to get an estimate down to the microscopic scale $O(\varepsilon)$. 

A key step in the proof of boundary Lipschitz estimates is to estimate boundary layer correctors, which is done by combining the classical Lipschitz estimate with a uniform H\"older estimate, as in \cite[Lemma 17]{alin} or \cite[Lemma 10]{BLtailosc}. We are able to relax the regularity assumption on $\psi$. This progress is enabled by our new estimate \eqref{estbumpyhpuloc} for the boundary layer corrector, which holds for Lipschitz boundaries $\psi$.

We begin with an estimate which is of constant use in this part of our work. Take $\psi\in W^{1,\infty}(\mathbb R^{d-1})$ and $A\in\mathcal A^{\nu}$. By Cacciopoli's inequality, there exists $C>0$ such that for all $\varepsilon>0$, for all weak solution $u^\varepsilon$ to
\begin{equation}\label{sysuepslemstep1}
\left\{\begin{array}{rll}
-\nabla\cdot A(x/\varepsilon)\nabla u^\varepsilon&=0,&x\in D^\varepsilon_\psi(0,1),\\
u^\varepsilon&=0,&x\in\Delta^\varepsilon_\psi(0,1),
\end{array}\right.
\end{equation}
for all $0<\theta<1$, 
\begin{equation}\label{estcacciomeanpartialdueps}
\begin{aligned}
\left|(\overline{\partial_{x_d} u^\varepsilon})_{D_\psi(0,\theta)}\right|&=\left|\intbar_{D_\psi(0,\theta)}\partial_{x_d}u^\varepsilon\right|\leq \left(\intbar_{D_\psi(0,\theta)}|\partial_{x_d}u^\varepsilon|^2\right)^{1/2}\\
&\leq \frac{C_0}{\theta^{d/2}(1-\theta)}\left(\intbar_{D_\psi(0,1)}|u^\varepsilon|^2\right)^{1/2}.
\end{aligned}
\end{equation}
Notice that $C_0$ in \eqref{estcacciomeanpartialdueps} only depends on $\lambda$.

Proposition \ref{proplipdowmicro} is a consequence of the two following lemmas. The first one contains the compactness argument. The second one is the iteration lemma. In order to alleviate the statement of the following lemma, the definition of the boundary layer $v$ is given straight after the lemma.

\begin{lem}\label{lemstep1lipdownmicro}
For all $\nu>0$, $\gamma>0$, there exists $\theta>0$, $0<\mu<1$, $\varepsilon_0>0$, such that for all $\psi\in W^{1,\infty}(\mathbb R^{d-1})$, $-1<\psi<0$ and $\|\nabla\psi\|_{L^\infty}\leq\gamma$, for all $A\in\mathcal A^{\nu}$, for all $0<\varepsilon<\varepsilon_0$, for all weak solution $u^\varepsilon$ to \eqref{sysuepslemstep1} we have
\begin{equation*}
\intbar_{D^\varepsilon_\psi(0,1)}|u^\varepsilon|^2\leq 1
\end{equation*}
implies
\begin{equation*}
\intbar_{D^\varepsilon_\psi(0,\theta)}\left|u^\varepsilon(x)-(\overline{\partial_{x_d} u^\varepsilon})_{D^\varepsilon_\psi(0,\theta)}\left[x_d+\varepsilon\chi^d(x/\varepsilon)+\varepsilon v(x/\varepsilon)\right]\right|^2dy\leq \theta^{2+2\mu}.
\end{equation*}
\end{lem}

The boundary layer $v=v(y)$ is the unique solution given by Theorem \ref{theoblbumpy} to the system
\begin{equation}\label{blsysiteration}
\left\{
\begin{array}{rll}
-\nabla\cdot A(y)\nabla v&=0,&y_d>\psi(y'),\\
v&=y_d+\chi^d(y),&y_d=\psi(y').
\end{array}
\right.
\end{equation}
The estimate of Theorem \ref{theoblbumpy} implies
\begin{equation*}
\sup_{\xi\in\mathbb Z^{d-1}}\int_{\xi+(0,1)^{d-1}}\int_{\psi(y')}^\infty|\nabla v|^2\leq C\left\{\|\psi\|_{H^{1/2}_{uloc}(\mathbb R^{d-1})}+\|\chi(\cdot,\psi(\cdot))\|_{H^{1/2}_{uloc}(\mathbb R^{d-1})}\right\},
\end{equation*}
with $C=C(d,N,\lambda,[A]_{C^{0,\nu}},\|\psi\|_{W^{1,\infty}})$. Now, by Sobolev injection $W^{1,\infty}(\mathbb R^{d-1})\hookrightarrow H^{1/2}(\mathbb R^{d-1})$ 
\begin{equation*}
\|\psi\|_{H^{1/2}_{uloc}(\mathbb R^{d-1})}\leq C\|\psi\|_{W^{1,\infty}(\mathbb R^{d-1})},
\end{equation*}
with $C=C(d)$ and by classical interior Lipschitz regularity
\begin{equation*}
\|\chi(\cdot,\psi(\cdot))\|_{H^{1/2}_{uloc}}\leq C\|\chi(\cdot,\psi(\cdot))\|_{W^{1,\infty}(\mathbb R^{d-1})}\leq C\|\chi\|_{W^{1,\infty}(\mathbb R^d)}\leq C,
\end{equation*}
with in the last inequality $C=C(d,N,\lambda,[A]_{C^{0,\nu}})$. Eventually,
\begin{equation}\label{mainestblcorrector}
\sup_{\xi\in\mathbb Z^{d-1}}\int_{\xi+(0,1)^{d-1}}\int_{\psi(y')}^\infty|\nabla v|^2\leq C,
\end{equation}
with $C=C(d,N,\lambda,[A]_{C^{0,\nu}},\|\psi\|_{W^{1,\infty}})$ uniform in $\varepsilon$.

\begin{lem}\label{lemstep2lipdownmicro}
Let $\theta$, $\varepsilon_0$ and $\gamma$ be given as in Lemma \ref{lemstep1lipdownmicro}. For all $\psi\in W^{1,\infty}(\mathbb R^{d-1})$, $-1<\psi<0$ and $\|\nabla\psi\|_{L^\infty}\leq\gamma$, for all $A\in\mathcal A^{\nu}$, for all $k\in\mathbb N$, $k>0$, for all $0<\varepsilon<\theta^{k-1}\varepsilon_0$, for all weak solution $u^\varepsilon$ to \eqref{sysuepslemstep1} there exists $a^\varepsilon_k\in\mathbb R^N$ satifying
\begin{equation*}
|a^\varepsilon_k|\leq C_0\frac{1+\theta^\mu+\ldots\ \theta^{\mu(k-1)}}{\theta^{d/2}(1-\theta)},
\end{equation*}
such that
\begin{equation*}
\intbar_{D^\varepsilon_\psi(0,1)}|u^\varepsilon|^2\leq 1
\end{equation*}
implies
\begin{equation}\label{estitstep2orderk}
\intbar_{D^\varepsilon_\psi(0,\theta^k)}\left|u^\varepsilon(x)-a^\varepsilon_k\left[x_d+\varepsilon\chi^d(x/\varepsilon)+\varepsilon v(x/\varepsilon)\right]\right|^2dy\leq \theta^{(2+2\mu)k},
\end{equation}
where $v=v(y)$ is the solution, given by Theorem \ref{theoblbumpy}, to the boundary layer system \eqref{blsysiteration}.
\end{lem}

The condition $\varepsilon<\theta^{k-1}\varepsilon_0$ can be seen as giving a lower bound on the scales $\theta^k$ for which one can prove the regularity estimate: $\theta^{k-1}>\varepsilon/\varepsilon_0$. In that perspective, estimate \eqref{estitstep2orderk} is an improved $C^{1,\mu}$ estimate down to the microscale $\varepsilon/\varepsilon_0$. 

For fixed $0<\varepsilon/\varepsilon_0<1/2$ and $r\in[\varepsilon/\varepsilon_0,1/2]$, there exists $k\in\mathbb N$ such that $\theta^{k+1}<r\leq\theta^k$. We aim at estimating 
\begin{equation*}
\intbar_{D^\varepsilon_\psi(0,r)}|u^\varepsilon(x)|^2
\end{equation*}
using the bound \eqref{estitstep2orderk}. We have
\begin{equation}\label{splittrianglefinalest}
\begin{aligned}
&\left(\intbar_{D^\varepsilon_\psi(0,r)}|u^\varepsilon(x)|^2\right)^{1/2}\leq\left(\intbar_{D^\varepsilon_\psi(0,\theta^k)}|u^\varepsilon(x)|^2\right)^{1/2}\\
&\leq\left(\intbar_{D^\varepsilon_\psi(0,\theta^k)}\left|u^\varepsilon(x)-a^\varepsilon_k\left[x_d-\psi(x')+\varepsilon\chi^d(x/\varepsilon)+\varepsilon v(x/\varepsilon)\right]\right|^2dy\right)^{1/2}\\
&\qquad\qquad+|a^\varepsilon_k|\left\{\left(\intbar_{D^\varepsilon_\psi(0,\theta^k)}|x_d|^2\right)^{1/2}+\left(\intbar_{D^\varepsilon_\psi(0,\theta^k)}|\varepsilon\chi^d(x/\varepsilon)|^2\right)^{1/2}+\left(\intbar_{D^\varepsilon_\psi(0,\theta^k)}|\varepsilon v(x/\varepsilon)|^2\right)^{1/2}\right\}.
\end{aligned}
\end{equation}
Let us focus on the term involving the boundary layer. Let $\eta=\eta(y_d)\in C^\infty_c(\mathbb R)$ be a cut-off such that $\eta\equiv 1$ on $(-1,1)$ and $\supp\eta\subset (-2,2)$. The triangle inequality yields
\begin{multline*}
\left(\intbar_{D^\varepsilon_\psi(0,\theta^k)}|\varepsilon v(x/\varepsilon)|^2\right)^{1/2}\leq \left(\intbar_{D^\varepsilon_\psi(0,\theta^k)}|\varepsilon v(x/\varepsilon)-(x_d+\varepsilon\chi^d(x/\varepsilon))\eta(x_d/\varepsilon)|^2\right)^{1/2}\\
+\left(\intbar_{D^\varepsilon_\psi(0,\theta^k)}|(x_d+\varepsilon\chi^d(x/\varepsilon))\eta(x_d/\varepsilon)|^2\right)^{1/2}.
\end{multline*}
Poincar\'e's inequality implies
\begin{equation*}
\begin{aligned}
&\left(\intbar_{D^\varepsilon_\psi(0,\theta^k)}|\varepsilon v(x/\varepsilon)-(x_d+\varepsilon\chi^d(x/\varepsilon))\eta(x_d/\varepsilon)|^2\right)^{1/2}\\
&\leq \theta^k\left(\intbar_{D^\varepsilon_\psi(0,\theta^k)}\left|\nabla\left(\varepsilon v(x/\varepsilon)-(x_d+\varepsilon\chi^d(x/\varepsilon))\eta(x_d/\varepsilon)\right)\right|^2\right)^{1/2}\\
&\leq \theta^k\left(\intbar_{D^\varepsilon_\psi(0,\theta^k)}|\nabla v(x/\varepsilon)|^2\right)^{1/2}+(1+\|\nabla\chi\|_{L^\infty}^2)\theta^k\left(\intbar_{D^\varepsilon_\psi(0,\theta^k)}|\eta(x_d/\varepsilon)|^2\right)^{1/2}\\
&\qquad\qquad+\frac{\theta^k}{\varepsilon}\left(\intbar_{D^\varepsilon_\psi(0,\theta^k)}|(x_d+\varepsilon\chi^d(x/\varepsilon))\eta'(x_d/\varepsilon)|^2\right)^{1/2}.
\end{aligned}
\end{equation*}
Estimate \eqref{mainestblcorrector} now yields
\begin{equation*}
\intbar_{D^\varepsilon_\psi(0,\theta^k)}|\nabla v(x/\varepsilon)|^2\leq C\varepsilon\theta^{-k},
\end{equation*}
so that eventually using $\varepsilon/\varepsilon_0\leq r\leq\theta^k$,
\begin{equation*}
\left(\intbar_{D^\varepsilon_\psi(0,\theta^k)}|\varepsilon v(x/\varepsilon)-(x_d+\varepsilon\chi^d(x/\varepsilon))\eta(x_d/\varepsilon)|^2\right)^{1/2}\leq C\left(\varepsilon^{1/2}\theta^{k/2}+\theta^k+\frac{\theta^k}{\varepsilon}(\varepsilon+\varepsilon)\right)\leq C\theta^k
\end{equation*}
with $C=C(d,N,\lambda,[A]_{C^{0,\nu}},\|\psi\|_{W^{1,\infty}})$. It follows from \eqref{splittrianglefinalest} and \eqref{estitstep2orderk} that
\begin{equation*}
\left(\intbar_{D^\varepsilon_\psi(0,r)}|u^\varepsilon(x)|^2\right)^{1/2}\leq\theta^{(1+\mu)k}+C\theta^k\leq C\theta^k\leq Cr,
\end{equation*}
which is the estimate of Proposition \ref{proplipdowmicro}.

\subsection{Proof of Lemma \ref{lemstep1lipdownmicro}}

Let $0<\theta<1/8$ and $u^0\in H^1(D_0(0,1/4))$ be a weak solution of 
\begin{equation}\label{sysu0}
\left\{
\begin{array}{rll}
-\nabla\cdot A^0\nabla u^0&=0,&x\in D_0(0,1/4),\\
u^0&=0,&x\in \Delta_0(0,1/4),
\end{array}
\right. 
\end{equation}
such that 
\begin{equation*}
\intbar_{D_0(0,1/4)}|u^0|^2\leq 4^d. 
\end{equation*}
The classical regularity theory yields $u^0\in C^2(\overline{D_0(0,1/8)})$. Using that for all $x\in D_0(0,\theta)$
\begin{equation*}
\begin{aligned}
u^0(x)-\left(\overline{\partial_{x_d} u^0}\right)_{0,\theta}x_d&=u^0(x)-u^0(x',0)-\left(\overline{\partial_{x_d} u^0}\right)_{0,\theta}x_d\\
&=\frac{1}{|D^0(0,\theta)|}\int_0^1\int_{D_0(0,\theta)}\left(\partial_{x_d}u^0(x',tx_d)-\partial_{x_d}u^0(y)\right)x_ddydt.
\end{aligned}
\end{equation*}
we get
\begin{equation}\label{classSchauderu0}
\intbar_{D_0(0,\theta)}\left|u^0(x)-(\overline{\partial_{x_d} u^0})_{D_0(0,\theta)}x_d\right|^2dy\leq \widehat{C}\theta^4,
\end{equation}
where $\widehat{C}=\widehat{C}(d,N,\lambda)$. Fix $0<\mu<1$. Choose $0<\theta<1/8$ sufficiently small such that
\begin{equation}\label{estthetatocontra}
\theta^{2+2\mu}>\widehat{C}\theta^4.
\end{equation}

The rest of the proof is by contradiction. Fix $\gamma>0$. Assume that for all $k\in\mathbb N$, there exists $\psi_k\in W^{1,\infty}(\mathbb R^{d-1})$, 
\begin{equation}\label{condpsik}
-1<\psi<0\quad\mbox{and}\quad\|\psi_k\|_{L^\infty}\leq \gamma, 
\end{equation}
there exists $A_k\in\mathcal A^\nu$, there exists $0<\varepsilon_k<1/k$, there exists $u^{\varepsilon_k}_k$ solving 
\begin{equation*}
\left\{\begin{array}{rll}
-\nabla\cdot A_k(x/\varepsilon_k)\nabla u^{\varepsilon_k}_k&=0,&x\in D^{\varepsilon_k}_{\psi_k}(0,1),\\
u^{\varepsilon_k}_k&=0,&x\in\Delta^{\varepsilon_k}_{\psi_k}(0,1),
\end{array}\right.
\end{equation*}
such that
\begin{equation}\label{contrastep1hyp}
\intbar_{D^{\varepsilon_k}_{\psi_k}(0,1)}|u^{\varepsilon_k}_k|^2\leq 1
\end{equation}
and
\begin{equation}\label{contrastep1est}
\intbar_{D^{\varepsilon_k}_{\psi_k}(0,\theta)}\left|u^{\varepsilon_k}_k(x)-(\overline{\partial_{x_d} u^{\varepsilon_k}_k})_{D^{\varepsilon_k}_{\psi_k}(0,\theta)}\left[x_d+\varepsilon_k\chi^d_k(x/\varepsilon_k)+\varepsilon_k v_k(x/\varepsilon_k)\right]\right|^2dy>\theta^{2+2\mu}.
\end{equation}
Notice that $\chi^d_k$ is the cell corrector associated to the operator $-\nabla\cdot A_k(y)\nabla$ and $v_k$ is the boundary layer corrector associated to $-\nabla\cdot A_k(y)\nabla$ and to the domain $y_d>\psi_k(y')$.

First of all, for technical reasons, let us extend $u^{\varepsilon_k}_k$ by zero below the boundary, on $\{x'\in(-1,1)^{d-1},\ x_d\leq\varepsilon_k\psi_k(x'/\varepsilon_k)\}$. The extended functions are still denoted the same, and $u^{\varepsilon_k}_k$ is a weak solution of
\begin{equation*}
-\nabla\cdot A_k(x/\varepsilon_k)\nabla u^{\varepsilon_k}_k=0
\end{equation*}
on $\{x'\in(-1,1)^{d-1},\ x_d\leq\varepsilon_k\psi_k(x'/\varepsilon_k)+1\}$.

For $k$ sufficiently large, by Cacciopoli's inequality,
\begin{equation*}
\int_{(-1/4,1/4)^{d}}|\nabla u^{\varepsilon_k}_k|^2\leq C\int_{(-1/2,1/2)^{d}}|u^{\varepsilon_k}_k|^2\leq C,
\end{equation*}
where $C=C(d,N,\lambda)$. Therefore, up to a subsequence, which we denote again by $u^{\varepsilon_k}_k$, we have
\begin{equation}\label{cvstep1uepsk}
\begin{aligned}
u^{\varepsilon_k}_k\stackrel{k\rightarrow\infty}{\longrightarrow}u^0,\quad\mbox{strongly in}\quad L^2((-1/4,1/4)^{d-1}\times(-1,1/4)),\\
\nabla u^{\varepsilon_k}_k\stackrel{k\rightarrow\infty}{\rightharpoonup}\nabla u^0,\quad\mbox{weakly in}\quad L^2((-1/4,1/4)^{d-1}\times(-1,1/4)).
\end{aligned}
\end{equation}
Moreover, $\varepsilon_k\psi_k(\cdot/\varepsilon_k)$ converges to $0$ because $\psi_k$ is bounded uniformly in $k$ (see \eqref{condpsik}). Let $\varphi\in C^\infty_c(D_0(0,1/4))$. Theorem \ref{theoweakcvhomo} implies that 
\begin{equation*}
\int_{D^{\varepsilon_k}_{\psi_k}(0,1/4)}A_k(x/\varepsilon_k)\nabla u^{\varepsilon_k}_k\cdot\nabla\varphi\stackrel{k\rightarrow\infty}{\longrightarrow}\int_{D_0(0,1/4)}A^0\nabla u^0\nabla\varphi,
\end{equation*}
so that $u^0$ is a weak solution to
\begin{equation*}
-\nabla\cdot A^0\nabla u^0=0\quad\mbox{in}\quad D_0(0,1/4).
\end{equation*}
Furthermore, for all $\varphi\in C^\infty_c((-1/4,1/4)^{d-1}\times(-1,0))$, 
\begin{equation*}
0=\int_{\{x'\in(-1/4,1/4)^{d-1},\ -1\leq x_d\leq\varepsilon_k\psi_k(x'/\varepsilon_k)\}}u^{\varepsilon_k}_k\varphi\stackrel{k\rightarrow\infty}{\longrightarrow}\int_{(-1/4,1/4)^{d-1}\times(-1,0)}u^0\varphi,
\end{equation*}
so that $u^0(x)=0$ for all $x\in (-1/4,1/4)^{d-1}\times(-1,0)$. In particular, $u^0=0$ in $H^{1/2}(\Delta_0(0,1/4))$. Thus, $u^0$ is a solution to \eqref{sysu0} and satisfies the estimate \eqref{classSchauderu0}.

It remains to pass to the limit in \eqref{contrastep1est} to reach a contradiction. Since $|D^{\varepsilon_k}_{\psi_k}(0,\theta)|=|D_0(0,\theta)|$, we have
\begin{multline}\label{cvmeanstep1}
\left|(\overline{\partial_{x_d}u^{\varepsilon_k}_k})_{D^{\varepsilon_k}_{\psi_k}(0,\theta)}-(\overline{\partial_{x_d}u^0})_{D_0(0,\theta)}\right|\leq \frac{1}{|D_0(0,\theta)|}\left[\left|\int_{D^{\varepsilon_k}_{\psi_k}(0,\theta)\cap D_0(0,\theta)}\left(\partial_{x_d}u^{\varepsilon_k}_k-\partial_{x_d}u^0\right)\right|\right.\\
\left.+\int_{\left(D^{\varepsilon_k}_{\psi_k}(0,\theta)\setminus D_0(0,\theta)\right)\cup\left(D_0(0,\theta)\setminus D^{\varepsilon_k}_{\psi_k}(0,\theta)\right)}\left|\partial_{x_d}u^{\varepsilon_k}_k-\partial_{x_d}u^0\right|\right].
\end{multline}
The first term in the right hand side of \eqref{cvmeanstep1} tends to $0$ thanks to the weak convergence of $\nabla u^{\varepsilon_k}_k$ in \eqref{cvstep1uepsk}. The second term in the right hand side of \eqref{cvmeanstep1} goes to $0$ when $k\rightarrow\infty$ because of the $L^2$ bound on the gradient, and the fact that 
\begin{equation*} 
\left|\left(D^{\varepsilon_k}_{\psi_k}(0,\theta)\setminus D_0(0,\theta)\right)\cup\left(D_0(0,\theta)\setminus D^{\varepsilon_k}_{\psi_k}(0,\theta)\right)\right|\stackrel{k\rightarrow\infty}{\longrightarrow} 0.
\end{equation*}
Therefore,
\begin{equation*}
\intbar_{D^{\varepsilon_k}_{\psi_k}(0,\theta)\cap D_0(0,\theta)}\left|(\overline{\partial_{x_d} u^{\varepsilon_k}_k})_{D^{\varepsilon_k}_{\psi_k}(0,\theta)}\left[x_d+\varepsilon_k\chi^d_k(x/\varepsilon_k)\right]-(\overline{\partial_{x_d}u^0})_{D_0(0,\theta)}x_d\right|^2\stackrel{k\rightarrow\infty}{\longrightarrow}0.
\end{equation*}
Moreover, the strong $L^2$ convergence in \eqref{cvstep1uepsk} implies
\begin{equation*}
\intbar_{D^{\varepsilon_k}_{\psi_k}(0,\theta)\cap D_0(0,\theta)}|u^{\varepsilon_k}_k-u^0|^2\stackrel{k\rightarrow\infty}{\longrightarrow}0.
\end{equation*}
The last thing we have to check is the convergence 
\begin{equation*}
\intbar_{D^{\varepsilon_k}_{\psi_k}(0,\theta)}|\varepsilon_kv_k(x/\varepsilon_k)|^2\stackrel{k\rightarrow\infty}{\longrightarrow}0.
\end{equation*}
Let $\eta=\eta(y_d)\in C^\infty_c(\mathbb R)$ such that $\eta\equiv 1$ on $(-1,1)$ and $\supp\eta\subset(-2,2)$. We have
\begin{equation}\label{controlstep1termvk}
\begin{aligned}
&\intbar_{D^{\varepsilon_k}_{\psi_k}(0,\theta)}|\varepsilon_kv_k(x/\varepsilon_k)|^2\\
&\leq\intbar_{D^{\varepsilon_k}_{\psi_k}(0,\theta)}|\varepsilon_kv_k(x/\varepsilon_k)-(x_d+\varepsilon_k\chi^d_k(x/\varepsilon_k))\eta(x_d/\varepsilon_k)|^2+\intbar_{D^{\varepsilon_k}_{\psi_k}(0,\theta)}|(x_d+\varepsilon_k\chi^d_k(x/\varepsilon_k))\eta(x_d/\varepsilon_k)|^2.
\end{aligned}
\end{equation}
The last term in the right hand side of \eqref{controlstep1termvk} goes to $0$ when $k\rightarrow\infty$. Now by Poincar\'e's inequality,
\begin{equation*}
\begin{aligned}
&\intbar_{D^{\varepsilon_k}_{\psi_k}(0,\theta)}|\varepsilon_kv_k(x/\varepsilon_k)-(x_d+\varepsilon_k\chi^d_k(x/\varepsilon_k))\eta(x_d/\varepsilon_k)|^2\\
&\leq C\theta^2\left[\intbar_{D^{\varepsilon_k}_{\psi_k}(0,\theta)}|\nabla v_k(x/\varepsilon_k)|^2+\intbar_{D_{\psi_k}(0,\theta)}\left|\nabla\left((x_d+\varepsilon_k\chi^d_k(x/\varepsilon_k))\eta(x_d/\varepsilon_k)\right)\right|^2\right].
\end{aligned}
\end{equation*}
On the one hand by estimate \eqref{mainestblcorrector}
\begin{equation*}
\begin{aligned}
\intbar_{D^{\varepsilon_k}_{\psi_k}(0,\theta)}|\nabla v_k(x/\varepsilon_k)|^2&\leq \frac{C\varepsilon_k^d}{\theta^d}\int_{D^1_{\psi_k}(0,\theta/\varepsilon_k)}|\nabla v_k(y)|^2\\
&\leq C\varepsilon_k\sup_{\xi\in\mathbb Z^{d-1}}\int_{\xi+(0,1)^{d-1}}\int_{\psi_k(y')}^{\infty}|\nabla v_k|^2\leq C\varepsilon_k\stackrel{k\rightarrow\infty}{\longrightarrow}0
\end{aligned}
\end{equation*}
with in the last inequality $C=C(d,N,\lambda,[A]_{C^{0,\nu}})$ uniform in $\varepsilon$, and on the other hand
\begin{multline*}
\intbar_{D^{\varepsilon_k}_{\psi_k}(0,\theta)}\left|\nabla\left((x_d+\varepsilon_k\chi^d_k(x/\varepsilon_k))\eta(x_d/\varepsilon_k)\right)\right|^2\leq (1+\|\nabla\chi_k\|_{L^\infty}^2)\intbar_{D^{\varepsilon_k}_{\psi_k}(0,\theta)}|\eta(x_d/\varepsilon_k)|^2\\
+\frac{1}{\varepsilon_k^2}\intbar_{D^{\varepsilon_k}_{\psi_k}(0,\theta)}|(x_d+\varepsilon_k\chi^d_k(x/\varepsilon_k))\eta'(x_d/\varepsilon_k)|^2\leq C\varepsilon_k\stackrel{k\rightarrow\infty}{\longrightarrow}0,
\end{multline*}
with in the last inequality $C=C(d,N,\lambda,[A]_{C^{0,\nu}})$. These convergence results imply that passing to the limit in \eqref{contrastep1est} we get
\begin{multline*}
\theta^{2+2\mu}\leq\intbar_{D^{\varepsilon_k}_{\psi_k}(0,\theta)}\left|u^{\varepsilon_k}_k(x)-(\overline{\partial_{x_d} u^{\varepsilon_k}_k})_{D^{\varepsilon_k}_{\psi_k}(0,\theta)}\left[x_d+\varepsilon_k\chi^d_k(x/\varepsilon_k)+\varepsilon_k v_k(x/\varepsilon_k)\right]\right|^2dy\\
\stackrel{k\rightarrow\infty}{\longrightarrow}\intbar_{D_0(0,\theta)}\left|u^0(x)-(\overline{\partial_{x_d} u^0})_{D_0(0,\theta)}x_d\right|^2dy\leq\widehat{C}\theta^4,
\end{multline*}
which contradicts \eqref{estthetatocontra}.

\subsection{Proof of Lemma \ref{lemstep2lipdownmicro}}

The proof is by induction on $k$. The result for $k=1$ is true because of Lemma \ref{lemstep1lipdownmicro}. Let $k\in\mathbb N$, $k\geq 1$. Assume that for all $\psi\in W^{1,\infty}(\mathbb R^{d-1})$ such that $-1<\psi<0$ and $\|\nabla\psi\|_{L^\infty}\leq\gamma$, for all $A\in\mathcal A^{\nu}$, for all $k\in\mathbb N$, $k>0$, for all $0<\varepsilon<\theta^{k-1}\varepsilon_0$, for all weak solution $u^\varepsilon$ to \eqref{sysuepslemstep1} there exists $a^\varepsilon_k\in\mathbb R^N$ satifying
\begin{equation*}
|a^\varepsilon_k|\leq C_0\frac{1+\theta^\mu+\ldots\ \theta^{\mu(k-1)}}{\theta^{d/2}(1-\theta)},
\end{equation*}
such that
\begin{equation*}
\intbar_{D^\varepsilon_\psi(0,1)}|u^\varepsilon|^2\leq 1
\end{equation*}
implies
\begin{equation}\label{estitstep2orderkproof}
\intbar_{D^\varepsilon_\psi(0,\theta^k)}\left|u^\varepsilon(x)-a^\varepsilon_k\left[x_d+\varepsilon\chi^d(x/\varepsilon)+\varepsilon v(x/\varepsilon)\right]\right|^2dy\leq \theta^{(2+2\mu)k}.
\end{equation}
This is our induction hypothesis. 

Given $\psi\in W^{1,\infty}(\mathbb R^{d-1})$, $-1<\psi<0$ and $\|\nabla\psi\|_{L^\infty}\leq\gamma$ and $A\in\mathcal A^{\nu}$, $0<\varepsilon<\theta^{k-1}\varepsilon_0$ and a solution $u^\varepsilon$ to \eqref{sysuepslemstep1} such that
\begin{equation*}
\intbar_{D^\varepsilon_\psi(0,1)}|u^\varepsilon|^2\leq 1
\end{equation*}
we define 
\begin{equation*}
U^\varepsilon(x):=\frac{1}{\theta^{(1+\mu)k}}\left\{u^\varepsilon(\theta^kx)-a^\varepsilon_k\left[\theta^kx_d+\varepsilon\chi^d(\theta^kx/\varepsilon)+\varepsilon v(\theta^kx/\varepsilon)\right]\right\}
\end{equation*}
for all $x\in D^{\varepsilon/\theta^k}_{\psi}(0,1)$. The goal is to apply the estimate of Lemma \ref{lemstep2lipdownmicro} to $U^\varepsilon$. By the induction estimate \eqref{estitstep2orderkproof}, we have
\begin{equation*}
\intbar_{D^{\varepsilon/\theta^k}_{\psi}(0,1)}|U^\varepsilon|^2\leq 1.
\end{equation*}
Moreover, $U^\varepsilon$ solves the system
\begin{equation}\label{sysUepsstep2}
\left\{
\begin{array}{rll}
-\nabla\cdot A(\theta^kx/\varepsilon)\nabla U^\varepsilon&=0,&x\in D^{\varepsilon/\theta^k}_{\psi}(0,1),\\
U^\varepsilon&=0,&x\in \Delta^{\varepsilon/\theta^k}_{\psi}(0,1).
\end{array}
\right.
\end{equation}
The boundary layer $v$ solving \eqref{blsysiteration} has been designed for $U^\varepsilon$ to solve \eqref{sysUepsstep2}. It follows that $U^\varepsilon$ satisfies the assumptions of Lemma \ref{lemstep1lipdownmicro}. Therefore, for all $\varepsilon/\theta^k<\varepsilon_0$, we have
\begin{equation*}
\intbar_{D^{\varepsilon/\theta^k}_\psi(0,\theta)}\left|U^\varepsilon(x)-(\overline{\partial_{x_d} U^\varepsilon})_{D^{\varepsilon/\theta^k}_\psi(0,\theta)}\left[x_d+\frac{\varepsilon}{\theta^k}\chi^d(\theta^kx/\varepsilon)+\frac{\varepsilon}{\theta^k} v(\theta^kx/\varepsilon)\right]\right|^2dy\leq \theta^{2+2\mu}.
\end{equation*}
Eventually,
\begin{equation*}
\intbar_{D^\varepsilon_\psi(0,\theta^{k+1})}\left|u^\varepsilon(x)-a^\varepsilon_{k+1}\left[x_d+\varepsilon\chi^d(x/\varepsilon)+\varepsilon v(x/\varepsilon)\right]\right|^2dy\leq \theta^{(2+2\mu)(k+1)},
\end{equation*}
with
\begin{equation*}
a^\varepsilon_{k+1}:=a^\varepsilon_k+\theta^{\mu k}(\overline{\partial_{x_d} U^\varepsilon})_{D^{\varepsilon/\theta^k}_\psi(0,\theta)}
\end{equation*}
satisfying the estimate
\begin{equation*}
|a^\varepsilon_{k+1}|\leq C_0\frac{1+\theta^\mu+\ldots\ \theta^{\mu(k-1)}}{\theta^{d/2}(1-\theta)}+C_0\frac{\theta^{\mu k}}{\theta^{d/2}(1-\theta)}\leq C_0\frac{1+\theta^\mu+\ldots\ \theta^{\mu k}}{\theta^{d/2}(1-\theta)}.
\end{equation*}
This concludes the iteration step and proves Lemma \ref{lemstep2lipdownmicro}.


\bibliographystyle{alpha}

\bibliography{Lp.bib}

\begin{thebibliography}{GVM10}

\bibitem[ABZ13]{ABZ13_Kato}
T.~{Alazard}, N.~{Burq}, and C.~{Zuily}.
\newblock {Cauchy theory for the gravity water waves system with non localized
  initial data}.
\newblock {\em ArXiv e-prints}, May 2013.

\bibitem[AL87a]{alin}
M.~Avellaneda and F.-H. Lin.
\newblock Compactness methods in the theory of homogenization.
\newblock {\em Comm. Pure Appl. Math}, 40(6):803--847, 1987.

\bibitem[AL87b]{alinscal}
M.~Avellaneda and F.-H. Lin.
\newblock Homogenization of elliptic problems with {$L^p$} boundary data.
\newblock {\em Appl. Math. Optim.}, 15(2):93--107, 1987.

\bibitem[AL89a]{alin2}
M.~Avellaneda and F.-H. Lin.
\newblock Compactness methods in the theory of homogenization. {II}.
  {E}quations in nondivergence form.
\newblock {\em Comm. Pure Appl. Math.}, 42(2):139--172, 1989.

\bibitem[AL89b]{Alin90P}
M.~Avellaneda and F.-H. Lin.
\newblock Homogenization of {P}oisson's kernel and applications to boundary
  control.
\newblock {\em J. Math. Pures Appl. (9)}, 68(1):1--29, 1989.

\bibitem[AL91]{alinLp}
M.~Avellaneda and F.-. Lin.
\newblock {$L^p$} bounds on singular integrals in homogenization.
\newblock {\em Comm. Pure Appl. Math.}, 44(8-9):897--910, 1991.

\bibitem[AS14a]{ArmstrongShen14}
S.~N. {Armstrong} and Z.~{Shen}.
\newblock {Lipschitz estimates in almost-periodic homogenization}.
\newblock {\em ArXiv e-prints}, September 2014.

\bibitem[AS14b]{ArSmart14}
S.~N. {Armstrong} and C.~K. {Smart}.
\newblock {Quantitative stochastic homogenization of convex integral
  functionals}.
\newblock {\em ArXiv e-prints}, June 2014.

\bibitem[DGV11]{DaGV11}
A.-L. Dalibard and D.~G{\'e}rard-Varet.
\newblock Effective boundary condition at a rough surface starting from a slip
  condition.
\newblock {\em J. Differential Equations}, 251(12):3450--3487, 2011.

\bibitem[DP14]{DP_SC}
A.-L. Dalibard and C.~Prange.
\newblock Well-posedness of the {S}tokes-{C}oriolis system in the half-space
  over a rough surface.
\newblock {\em Anal. PDE}, 7(6):1253--1315, 2014.

\bibitem[GNO14]{GloriaNeukammOtto14}
A.~{Gloria}, S.~{Neukamm}, and F.~{Otto}.
\newblock {A regularity theory for random elliptic operators}.
\newblock {\em ArXiv e-prints}, September 2014.

\bibitem[GVM10]{DGVNMnoslip}
D.~G{\'e}rard-Varet and N.~Masmoudi.
\newblock Relevance of the slip condition for fluid flows near an irregular
  boundary.
\newblock {\em Comm. Math. Phys.}, 295(1):99--137, 2010.

\bibitem[KLS13]{KLSNeumann}
C.~E. Kenig, F.-H. Lin, and Z.~Shen.
\newblock Homogenization of elliptic systems with {N}eumann boundary
  conditions.
\newblock {\em J. Amer. Math. Soc.}, 26(4):901--937, 2013.

\bibitem[KP15]{BLtailosc}
C.~{Kenig} and C.~{Prange}.
\newblock {Uniform Lipschitz Estimates in Bumpy Half-Spaces}.
\newblock {\em Arch. Rational Mech. Anal.}, 216(3):703--765, 2015.

\bibitem[LS80]{LS}
O.~A. Lady{\v{z}}enskaja and V.~A. Solonnikov.
\newblock Determination of solutions of boundary value problems for stationary
  {S}tokes and {N}avier-{S}tokes equations having an unbounded {D}irichlet
  integral.
\newblock {\em Zap. Nauchn. Sem. Leningrad. Otdel. Mat. Inst. Steklov. (LOMI)},
  96:117--160, 308, 1980.
\newblock Boundary value problems of mathematical physics and related questions
  in the theory of functions, 12.

\end{thebibliography}

\end{document}